\documentclass[
final,
11pt, a4paper, USenglish
]{article}

\usepackage{import}
\usepackage[USenglish]{babel}
\usepackage{import}

\usepackage[utf8]{inputenc}

\usepackage[T1]{fontenc}
\usepackage{cite}

\usepackage{mathscinet}

\usepackage{amsmath, amscd, amsthm}
\usepackage{thmtools}
\usepackage{amssymb}
\usepackage{mathtools}

\usepackage[shortcuts]{extdash}

\usepackage{tikz-cd}
\usepackage{tikz}

\usepackage{ifthen}

\newcommand{\ifisemptythenelse}[3]{
  \if\relax\detokenize{#1}\relax
    #2%
  \else
    #3%
  \fi
}
\newcommand{\maybebrackets}[1]{
				\ifisemptythenelse{#1}{}{{(#1)}}
				}
\newcommand{\maybesubscript}[2]{
  \ifisemptythenelse{#2}
  {#1}
  {{#1}_{#2}}
}

\usepackage{dsfont}
\usepackage[mathcal]{eucal}

\newcommand\Eucal{\mathcal}
\newcommand\ol{\overline}
\newcommand\ul{\underline}

\newcommand{\mathblackboardboldnew}{\mathds}

\newcommand{\BB}{\mathblackboardboldnew B}

\newcommand{\BN}{\mathblackboardboldnew N}

\newcommand{\BS}{\mathblackboardboldnew S}

\newcommand{\BZ}{\mathblackboardboldnew Z}

\newcommand{\newmathcal}{\CMcal}

\newcommand{\calF}{\newmathcal{F}}

\newcommand{\rmK}{\mathrm{K}}

\newcommand{\rmS}{\mathrm{S}}

\newcommand{\frakS}{\mathfrak{S}}

\newcommand{\qquote}[1]{``#1''}
\newcommand{\roughly}[1]{``#1''}

\newcommand{\introduce}{\textbf}
\newcommand{\latin}{\emph}
\newcommand{\buzzword}{\emph}

\makeatletter
\newcommand\defpureword[2]{%
  \expandafter\newcommand\csname #1\endcsname[1]{\@pureword{#2}{##1}{#1}}
}
\newcommand{\@pureword}[3]{
  \ifisemptythenelse%
  {#2}%
  {#1}%
  {
    \GenericError{}
    {Argument of pureword command not empty}
    {Use \ #3{} instead of \ #3.}
  }%
}
\newcommand\pureword[2]{\@pureword{#1}{#2}{command}}
\makeatother

\defpureword{Ktheory}{$\rmK$-theory}

\defpureword{Sconstruction}{$\rmS$-construction}

\defpureword{ConnesCC}{Connes' cyclic category}

\defpureword{inftycat}{$\infty$\=/category}
\defpureword{inftycats}{$\infty$\=/categories}
\defpureword{inftycategorical}{$\infty$\=/categorical}
\defpureword{pinftycat}{($\infty$\=/)category}
\defpureword{pinftycats}{($\infty$\=/)categories}
\defpureword{inftyoperad}{$\infty$\=/operad}
\defpureword{inftyoperads}{$\infty$\=/operads}
\defpureword{pinftyoperad}{{($\infty$\=/)}operad}
\defpureword{pinftyoperads}{{($\infty$\=/)}operads}

\newcommand{\resp}{{resp.\ }}
\newcommand\auxiliarytoggleablelatin\latin

\newcommand{\ie}{\auxiliarytoggleablelatin{i.e.}, }
\newcommand{\cf}{cf.\ }
\newcommand{\eg}{\auxiliarytoggleablelatin{e.g.}, }

\newcommand{\aka}{a.k.a.\ }

\newcommand{\PhD}{Ph.D.\ }

\newcommand{\non}[1]{non-{#1}}

\newcommand{\pre}[1]{pre{#1}}
\newcommand{\ppre}[1]{(pre){#1}}
\newcommand{\sub}[1]{sub{#1}}
\newcommand{\self}[1]{self-{#1}}
\newcommand{\well}[1]{well {#1}}

\newcommand{\twonames}[2]{{#1} and {#2}}

\defpureword{Cech}{\v{C}ech}

\defpureword{TDyckerhoff}{T.~Dyckerhoff}
\defpureword{GJasso}{G.~Jasso}

\defpureword{theauthor}{the author} 
\defpureword{Theauthor}{The author} 

\defpureword{htpyterminal}{homotopy terminal}
\defpureword{htpyinitial}{homotopy initial}
\defpureword{wcontractible}{weakly contractible}

\newtheoremstyle{newtheorem}{}{}{}{}{\bfseries}

\DeclareSymbolFont{extraup}{U}{zavm}{m}{n}
\DeclareMathSymbol{\varheart}{\mathalpha}{extraup}{86}
\DeclareMathSymbol{\vardiamond}{\mathalpha}{extraup}{87}
\DeclareMathSymbol{\varclub}{\mathalpha}{extraup}{84}
\DeclareMathSymbol{\varspade}{\mathalpha}{extraup}{85}

\newcommand{\somberend}{\ensuremath{\Diamond}}

\ifx\fancyends\undefined
\newcommand{\theoremend}{\ensuremath{\square}}
\newcommand{\definitionend}{\somberend}
\newcommand{\exampleend}{\somberend}
\newcommand{\questionend}{\somberend}
\newcommand{\proofend}{\ensuremath{\blacksquare}}
\else
\newcommand{\theoremend}{\ensuremath{\varheart}}
\newcommand{\definitionend}{\ensuremath{\clubsuit}}
\newcommand{\exampleend}{\ensuremath{\diamondsuit}}
\newcommand{\questionend}{\ensuremath{\spadesuit}}
\newcommand{\proofend}{\ensuremath{\square}}
\fi

\theoremstyle{definition}

\newtheorem{Masterthm}{Masterthm}[subsection]

\declaretheorem[name=Construction ,style=definition,qed={\definitionend},sibling=Masterthm]{Cstr}
\declaretheorem[name=Construction ,style=definition,qed={\definitionend},unnumbered]{Cstr*}
\declaretheorem[name=Corollary,style=definition,qed={\proofend},sibling=Masterthm]{dCor}
\declaretheorem[name=Corollary,style=definition,qed={\theoremend},sibling=Masterthm]{Cor}
\declaretheorem[name=Definition,style=definition,qed={\definitionend},sibling=Masterthm]{Def}
\declaretheorem[name=Definition,style=definition,qed={\definitionend},unnumbered]{Def*}

\declaretheorem[name=Goal,style=definition,qed={\questionend},unnumbered]{Goal*}

\declaretheorem[name=Lemma,style=definition,qed={\theoremend},sibling=Masterthm]{Lem}
\declaretheorem[name=Notation,style=definition,qed={\definitionend},sibling=Masterthm]{Not}

\declaretheorem[name=Proof,style=definition,qed={\proofend}, numbered=no]{Prf}

\declaretheorem[name=Proposition,style=definition,qed={\theoremend},sibling=Masterthm]{Prop}
\declaretheorem[name=Question,style=definition,qed={\questionend},sibling=Masterthm]{Qstn}

\declaretheorem[name=Theorem,style=definition,qed={\theoremend},sibling=Masterthm]{Thm}

\declaretheorem[name=Theorem,style=definition,qed={\theoremend},unnumbered]{Thm*}

\theoremstyle{remark}
\declaretheorem[name=Remark ,style=remark,qed={\exampleend},sibling=Masterthm]{Rem}

\declaretheorem[name=Fact,style=remark,qed={\exampleend},unnumbered]{Fact*}

\declaretheorem[name=Example ,style=remark,qed={\exampleend},sibling=Masterthm]{Expl}


\numberwithin{equation}{section}

\declaretheorem[name=Theorem, style=definition, qed={\somberend}]{Theorem}
\declaretheorem[name=Theorem, style=definition, qed={\somberend},unnumbered]{Theorem*}

\declaretheorem[name=Corollary, style=definition, qed={\somberend},unnumbered]{Corollary*}

\declaretheorem[name=Observation, style=definition, qed={\somberend}, unnumbered]{Observation*}

\declaretheorem[name=Conjecture, style=definition, qed={\somberend},unnumbered]{Conjecture*}

\declaretheorem[name=Definition, style=definition, qed={\somberend}, unnumbered]{Definition*}

\declaretheorem[name=Remark ,style=remark,qed={\somberend},unnumbered]{Remark*}

\ifx\nopagestyle\undefined
\usepackage{geometry}
\usepackage{enumitem}

\setlist{itemsep=-3pt, topsep=1pt}
\setlist[enumerate, 1]{label=(\arabic*), ref=(\arabic*)}

\newcommand{\axiomlabelstyle}[1]{(\bfseries{#1}\arabic*)} 
\newcommand{\axiomrefstyle}[1]{(#1\arabic*)} 

\usepackage{footnote}

\def\defaultdepthtoc{3}

\usepackage[titletoc]{appendix}				%
\usepackage[nottoc,notlot,notlof]{tocbibind} 	
\addto\captionsUSenglish{
  \renewcommand{\contentsname}%
    {Contents:}%
}
\addtocontents{toc}{\protect\setcounter{tocdepth}{\defaultdepthtoc}}

\makeatletter
\newcommand*\mytitle[1]{\gdef\@mytitle{#1}}
\newcommand*\myauthor[1]{\gdef\@myauthor{#1}}
\newcommand*\myaffiliation[1]{\gdef\@myaffiliation{#1}}
\newcommand*\myurl[1]{\gdef\@myurl{#1}}

\newcommand*\msnpError[1]{
  \GenericError{}
  {no my#1 provided}
  {Please call my#1 with an argument; do not try to def my#1}
}

\mytitle{\msnpError{title}}
\myauthor{\msnpError{author}}
\myaffiliation{\msnpError{affiliation}}
\myurl{\msnpError{url}}

\newcommand*\makemytitle{
  \title{\@mytitle}
  \author{\@myauthor\thanks{\@myaffiliation, \url{\@myurl}}}
  \maketitle
}

\usepackage{fancyhdr}
\usepackage{lastpage}

\makeatletter

\fi

\newcommand{\la}{\leftarrow}
\newcommand{\ra}{\rightarrow}

\newcommand{\hla}{\hookleftarrow}
\newcommand{\hra}{\hookrightarrow}
\newcommand{\thra}
{
  \mathrel
  {
    \mathchoice
    {\mathrlap{\rightarrow}\mkern -2mu\rightarrow}
    {\mathrlap{\rightarrow}\mkern -2mu\rightarrow}
    {\mathrlap{\rightarrow}\mkern 4mu\rightarrow}
    {\mathrlap{\rightarrow}\mkern 4mu\rightarrow}
  }
}

\newcommand{\lra}{\longrightarrow}
\newcommand{\lla}{\longleftarrow}
\newcommand{\llra}{\longleftrightarrow}

\newcommand{\lhla}{\longleftarrow\joinrel\rhook}
\newcommand{\lhra}{\lhook\joinrel\longrightarrow}
\newcommand{\lthra}
{
  \mathrel
  {
    \mathchoice
    {\mathrlap{\longrightarrow}\mkern 9mu\rightarrow}
    {\mathrlap{\longrightarrow}\mkern 9mu\rightarrow}
    {\mathrlap{\longrightarrow}\mkern 15mu\rightarrow}
    {\mathrlap{\longrightarrow}\mkern 15mu\rightarrow}
  }
}

\makeatletter

\newlength{\@minLengthArrowLengthOne}
\newlength{\@minLengthArrowLengthTwo}
\newlength{\@subscriptLengthForLongArrow}
\newlength{\@calculatedArrowLength}

\newcommand\@minLengthArrow [4][]{
  {
    \setlength{\@subscriptLengthForLongArrow}{#4}%
    \settowidth{\@minLengthArrowLengthOne}{\scriptsize$#1$}
    \settowidth{\@minLengthArrowLengthTwo}{\scriptsize$#2$}
    \pgfmathsetlength{\@calculatedArrowLength}{
      max
      ( \@minLengthArrowLengthOne
      , \@minLengthArrowLengthTwo
      , \@subscriptLengthForLongArrow
      )
    }
    \mathrel{
      #3
    [{\mathmakebox[\@calculatedArrowLength]{#1}}]
    {\mathmakebox[\@calculatedArrowLength]{#2}}
}  }
}

\newcommand\convertToMinLengthArrow[3][11pt]{
  \expandafter\newcommand\csname #2\endcsname[2][]{%
    \mathchoice
    {
      \mathrel{\@minLengthArrow[##1]{##2}{#3}{#1}}
    }
    {
      \mathrel{#3[##1]{##2}}
    }
    {
      \mathrel{#3[##1]{##2}}
    }
    {
      \mathrel{#3[##1]{##2}}
    }
  }
}

\makeatother

\convertToMinLengthArrow{xRa}{\xRightarrow}
\convertToMinLengthArrow{xLa}{\xLeftarrow}
\convertToMinLengthArrow{xLRa}{\xLeftrightarrow}
\convertToMinLengthArrow{xhla}{\xhookleftarrow}
\convertToMinLengthArrow{xhra}{\xhookrightarrow}
\convertToMinLengthArrow{xra}{\xrightarrow}
\convertToMinLengthArrow{xla}{\xleftarrow}
\convertToMinLengthArrow[13pt]{xlra}{\xleftrightarrow}
\convertToMinLengthArrow{xrat}{\xrightarrowtail}

\makeatletter
\newbox\xrat@below
\newbox\xrat@above
\newcommand{\xrightarrowtail}[2][]{%
  \setbox\xrat@below=\hbox{\ensuremath{\scriptstyle #1}}%
  \setbox\xrat@above=\hbox{\ensuremath{\scriptstyle #2}}%
  \pgfmathsetlengthmacro{\xrat@len}{max(\wd\xrat@below,\wd\xrat@above)+.6em}%
  \mathrel{\tikz [>->,baseline=-.75ex]
                 \draw (0,0) -- node[below=-2pt] {\box\xrat@below}
                                node[above=-2pt] {\box\xrat@above}
                       (\xrat@len,0) ;}}
\makeatother

\newcommand{\lrlas}{
  \mathrel{\substack{\lra \\[-.65ex] \lla}}
}

\newcommand\adjarrows{\lrlas}
\newcommand\ladjarrows\adjarrows
\newcommand\radjarrows\adjarrows

\usepackage{xstring}

\newcommand\displayset[1]{\left\{\text{#1}\right\}}
\newcommand\intxt[1]{\qquad\text{#1}\qquad}

\newcommand{\tikzcdadjunction}[3][r]{
  \ar[#1,shift left=2.5,"#2"]
  \ar[#1,phantom,description,"\tiny{\bot}"]
  \ar[from=#1,shift left=2.5,"#3"]
}

\newcommand\cdadjunctionOpt[5]{
  \IfEqCase{#1}{
    {u}{
      \ar[#1,shift left=2.5,"#2",/utils/exec={\pgfkeysalso{#4}}]
      \ar[#1,phantom,description,"\tiny{\dashv}"]
      \ar[from=#1,shift left=2.5,"#3",/utils/exec={\pgfkeysalso{#5}}]
    }
    {d}{
      \ar[#1,shift right=2.5,"#2"',/utils/exec={\pgfkeysalso{#4}}]
      \ar[#1,phantom,description,"\tiny{\dashv}"]
      \ar[from=#1,shift right=2.5,"#3"',/utils/exec={\pgfkeysalso{#5}}]
    }
    {r}{
      \ar[#1,shift left=2.5,"#2",/utils/exec={\pgfkeysalso{#4}}]
      \ar[#1,phantom,description,"\tiny{\bot}"]
      \ar[from=#1,shift left=2.5,"#3",/utils/exec={\pgfkeysalso{#5}}]
    }
    {l}{
      \ar[#1,shift right=2.5,"#2"',/utils/exec={\pgfkeysalso{#4}}]
      \ar[#1,phantom,description,"\tiny{\bot}"]
      \ar[from=#1,shift right=2.5,"#3"',/utils/exec={\pgfkeysalso{#5}}]
    }
  }
}

\newcommand\isCartesian[1][dr]{\ar[{#1}, phantom, description, very near start,"\lrcorner"]}
\newcommand\iscoCartesian[1][dr]{\ar[{#1}, phantom, description, very near end,"\ulcorner"]}
\newcommand\isbiCartesian[1][dr]{\ar[{#1},phantom, description, "\square"]}

\newcommand\cdtriangleoverOpt[7][]{%
  \def\obja{#2}%
  \def\objb{#3}%
  \def\objc{#4}%
  \def\mora{\pgfkeysalso{#5}}
  \def\morb{\pgfkeysalso{#6}}
  \def\morc{\pgfkeysalso{#7}}
  \begin{tikzcd}[ampersand replacement=\&]
    \obja\ar[rr,/utils/exec=\mora] \& \& \objb\\
    \& \objc\ar[from=ul,/utils/exec=\morc]\ar[from=ur,/utils/exec=\morb] \&
  \end{tikzcd}%
}

\newcommand\cdtriangleOpt[7][]{%
  \def\obja{#2}%
  \def\objb{#3}%
  \def\objc{#4}%
  \def\mora{\pgfkeysalso{#5}}
  \def\morb{\pgfkeysalso{#6}}
  \def\morc{\pgfkeysalso{#7}}
  \begin{tikzcd}[ampersand replacement=\&]
    \& \objb\ar[dr,/utils/exec=\morb] \&\\
    \obja\ar[rr,/utils/exec=\morc]\ar[ru,/utils/exec=\mora] \& \& \objc
  \end{tikzcd}%
}

\newcommand{\cdtriangle}[7][]{%
  \cdtriangleOpt[#1]{#2}{#3}{#4}{"{#5}"}{"{#6}"}{"{#7}"'}
}

\newcommand\cdsquare[9][]{
  \cdsquareOpt[#1]{#2}{#3}{#4}{#5}{"#6"}{"#7"'}{"#8"}{"#9"'}
}

\newcommand\cdsquareNA[5][]{\cdsquare[#1]{#2}{#3}{#4}{#5}{}{}{}{}}

\newcommand\cdsquareOpt[9][]{%
  \def\obja{#2}%
  \def\objb{#3}%
  \def\objc{#4}%
  \def\objd{#5}%
  \def\mora{\pgfkeysalso{#6}}
  \def\morb{\pgfkeysalso{#7}}
  \def\morc{\pgfkeysalso{#8}}
  \def\mord{\pgfkeysalso{#9}}
  \begin{tikzcd}[ampersand replacement=\&]
    \obja\ar[r,/utils/exec=\mora]\ar[d,/utils/exec=\morb]%
    \IfEqCase{#1}{
      {C}{\isCartesian}
      {cC}{\iscoCartesian}
      {bC}{\isbiCartesian}
    }
    \&%
    \objb\ar[d,/utils/exec=\morc]%
    \\%
    \objc\ar[r,/utils/exec=\mord]%
    \&%
    \objd%
  \end{tikzcd}%
}

\DeclareMathOperator*{\colim}{colim}

\DeclareMathOperator{\Fun}{Fun}
\DeclareMathOperator{\Map}{Map}

\DeclareMathOperator{\RKE}{RKE}
\DeclareMathOperator{\LKE}{LKE}
\DeclareMathOperator{\Res}{Res}

\DeclareMathOperator{\Ker}{Ker}

\DeclareMathOperator{\cof}{cof}
\DeclareMathOperator{\totfib}{tot-fib}
\DeclareMathOperator{\totcof}{tot-cof}

\renewcommand{\Im}{\operatorname{Im}}
\DeclareMathOperator{\Ob}{Ob}

\DeclareMathOperator{\Ar}{Ar}

\newcommand\Id{\mathrm{Id}}
\newcommand\pt{\mathrm{pt}}

\newcommand\blank{-}

\newcommand\inv{{-1}}
\newcommand\dual{{\vee}}
\newcommand\op{\mathrm{op}}

\newcommand\basepoint\star
\newcommand\card[1]{\left|#1\right|}
\newcommand\disjunion{\mathbin{\dot\cup}}
\newcommand\emptytuple\varnothing
\newcommand\invimage[2]{{#1}^{\inv}#2}
\newcommand\preimage[2]{\invimage{#1}{\set{#2}}}
\newcommand\setmid{\,\middle|\,}
\newcommand\setP[2]{\left\{#1\setmid #2\right\}}
\newcommand\set[1]{\left\{#1\right\}}
\renewcommand\emptyset\varnothing

\newcommand\invs[1]{{#1}^\inv}

\newcommand\SG[1]{{\frakS_{#1}}}
\newcommand\spherespectrum{\BS}

\DeclarePairedDelimiter\parentheses{\lparen}{\rparen}
\newcommand\apair[2]{\langle\ifisemptythenelse{#1}\blank{#1};\ifisemptythenelse{#2}\blank{#2}\rangle}
\newcommand\spanby[1]{\langle{#1}\rangle}

\newcommand\pushout[1]{\sqcup_{#1}}
\newcommand\fiberproduct[1]{\times_{#1}}
\newcommand\terminal{\star}
\newcommand\initial{\varnothing}
\newcommand{\loops}{\Omega}
\newcommand\ladjto\dashv
\newcommand\radjto\vdash
\newcommand\overcat[3][]{{{#2}_{/#3}}}
\newcommand\undercat[3][]{{{#3}_{#2/}}}
\newcommand\localizedat[2]{#1[\invs{#2}]}
\newcommand\htpycat[1]{\mathrm{h}{#1}}
\newcommand\grpdcore[1]{{#1}^{\simeq}}
\newcommand\ptdobjof[1]{{#1}_\star}

\newcommand\lcone[1]{{#1}^\triangleleft}
\newcommand\rcone[1]{{#1}^\triangleright}
\newcommand\lrcone[1]{{#1}^{\triangleleft\triangleright}}
\newcommand{\factorproj}[1]{\pi_{#1}}
\newcommand{\factorincl}[1]{\delta^{#1}}

 \newcommand\rquot[2]{
        \mathchoice
            {
                \text{\raise1ex\hbox{${#1}$}}\Big/\lower1ex\hbox{${#2}$}%
            }
            {
                #1\,/\,#2
            }
            {
                #1/#2
            }
            {
                #1/#2
            }
    }
\newcommand{\doubleslash}{\mathbin{
  \mathchoice{\Big/\mkern-10mu\Big/}
    {/\mkern-6mu/}
    {/\mkern-5mu/}
    {/\mkern-5mu/}}}

\newcommand{\qquot}[2]{
  \mathchoice
  {
    \text{\raise1ex\hbox{${#1}$}}\doubleslash\lower1ex\hbox{${#2}$}%
  }
  {{#1}\doubleslash{#2}}
  {{#1}\doubleslash{#2}}
  {{#1}\doubleslash{#2}}
}

\newcommand\lquot[2]{
        \mathchoice
            {
                \lower1ex\hbox{${#2}$}\Big\backslash\text{\raise1ex\hbox{${#1}$}}%
            }
            {
                #2\backslash#1
            }
            {
                #2\backslash#1
            }
            {
                #1\backslash#2
            }
    }
    
 \newcommand\lrquot[3]{
        \mathchoice
            {
                \lower1ex\hbox{${#2}$}\Big\backslash\text{\raise1ex\hbox{${#1}$}}\Big/\lower1ex\hbox{${#3}$}%
            }
            {
                #2\backslash#1/#3
            }
            {
                #2\backslash#1/#2
            }
            {
                #1\backslash#2/#2
            }
    }

\newcommand*\noloc{%
        \nobreak
        \mskip6mu plus1mu
        \mathpunct{}%
        \nonscript
        \mkern-\thinmuskip
        {:}%
        \mskip2mu
        \relax
}

\newcommand\restr[3][]{{
  \left.\kern-\nulldelimiterspace 
  {#2} 
  \vphantom{\big|} 
  \right|_{#3}^{#1} 
  }}

\newcommand{\category}[1]{\mathrm{\mathbf{#1}}}
\newcommand{\cat}\category

\newcommand{\rMod}[1]{\cat{Mod-}#1}
\newcommand{\lRep}[2]{{#1}\cat{-Rep}_{#2}}

\newcommand\finset{\cat{Fin}}

\newcommand\finsetp{\finset_\star}
\newcommand\finsetpop{\finset^\op_\star}
\newcommand\finbij{\grpdcore{\finset}}

\newcommand{\Set}{\cat{Set}}
\newcommand{\SetP}{\ptdobjof{\cat{Set}}}
\newcommand{\Ab}{\cat{Ab}}

\newcommand{\inftycatcal}[1]{\Eucal{#1}} 

\newcommand\Catinfty{\cat{Cat_\infty}}
\newcommand\Spaces{\Eucal{S}}
\newcommand\SpacesP{\ptdobjof{\Eucal{S}}}
\newcommand\Spectra{{\Eucal{S} p}}

\newcommand\numD[1]{{[#1]}}

\newcommand\Deltaleq[1]{\Delta_{\leq #1}}
\newcommand{\Dop}{{\Delta^{\op}}}
\newcommand{\Dopleq}[1]{{\Delta^{\op}_{\leq#1}}} 

\newcommand\nerve[1]{\mathrm N\maybebrackets{#1}}
\newcommand\geomrealization[1]{\card{\ifisemptythenelse{#1}{\blank}{#1}}}

\newcommand{\cofacemap}[1]{\mathrm{d}^{#1}}

\newcommand{\codegmap}[1]{\mathrm{s}^{#1}}

\newcommand\join{\mathbin\star}

\newcommand{\simplex}[1]{\Delta^{#1}}

\usepackage{hyperref}

\usepackage{ifdraft}
\usepackage[obeyFinal]{todonotes}
\AtEndDocument{\listoftodos}

\addto\extrasUSenglish{%
}

\usepackage{subfiles}
\mathtoolsset{showonlyrefs}

\mytitle{Homotopy coherent theorems of Dold--Kan type}
\myauthor{Tashi Walde}
\myaffiliation{Rheinische Friedrich-Wilhelms-Universität Bonn}
\myurl{https://www.math.uni-bonn.de/people/twalde/}

\begin{document}
\makemytitle

\begin{abstract}

We establish a large class of
homotopy coherent Morita-equivalences of Dold--Kan type
relating diagrams with values in any
weakly idempotent complete additive $\infty$-category;
the guiding example is an $\infty$-categorical Dold--Kan correspondence
between the $\infty$-categories
of simplicial objects and
connective coherent chain complexes.

Our results generalize many known
1-categorical equivalences
such as the classical Dold--Kan correspondence,
Pirashvili's Dold--Kan type theorem for abelian $\Gamma$-groups
and, more generally,
the combinatorial categorical equivalences of Lack and Street.

\end{abstract}

\tableofcontents

\newpage
\section{Introduction}

\newcommand{\FIsharp}{\mathrm{FI}\sharp}
\newcommand{\exring}{R}
The classical \DKcorr{}~\cite{Dold1958, Kan1958}
is a remarkable equivalence of categories
\begin{equation}
  \label{eq:intro-DK-classical}
  \Fun(\Dop,\exabcat) \xlra{\simeq} \connChof{\exabcat}
\end{equation}
between simplicial objects in $\exabcat$
and connective chain complexes in $\exabcat$,
where $\exabcat$ is the category of abelian groups
or, more generally, any abelian category~\cite{DoldPuppe1961}.
In the past decades,
many related equivalences have been constructed
\cite{Pirashvili2000
  ,Slominska2004
  ,Slominska2011
  ,Helmstutler2014
  ,CEF2015
  ,LackStreet2015
}
where the simplex category $\Delta$
is replaced by other categories
which are of similar \roughly{combinatorial nature}.

The goal of this article is to
simultaneously generalize these equivalences
in the homotopy coherent context of \inftycats{}.
To this end we study categories $\localB$
equipped with the structure $\excooltechdef$
of a so-called \buzzword{\cooltech{}}
(see \autoref{def:cooltech});
to each such \cooltech{} $\excooltech$
we associate a pointed category
$\localAp=\localAp(\excooltech)$
and prove the following
homotopy coherent correspondence of Dold--Kan type:

\begin{Theorem}[\autoref{cor:main-thm-without-pt}]
  \label{thm:intro-main}
  For each weakly idempotent complete\footnote{
    weakly idempotent complete = closed under direct complements
  } additive\footnote{
    additive = has direct sums and is enriched in abelian groups}
  \inftycat{} $\exwaddcat$,
  the \cooltech{} $\excooltech$
  induces a natural\footnote{
    natural in $\exwaddcat$ with respect to additive functors
  } equivalence
  \begin{equation}
    \label{eq:intro-main}
    \Fun(\localB,\exwaddcat)
    \xlra{\simeq}
    \ptdFun{\localAp}\exwaddcat
  \end{equation}
  between the \inftycats{} of diagrams $\localB\to\exwaddcat$
  and of pointed diagrams $\localAp\to\exwaddcat$.
\end{Theorem}

Before going into more details about \cooltechs{},
we explain how \autoref{thm:intro-main}
subsumes and generalizes previous results in the literature.

\begin{enumerate}
\item
  In the case where $\exabcat$ is an abelian category,
  we recover the classical \DKcorr{} \eqref{eq:intro-DK-classical}
  by applying \autoref{thm:intro-main} to $\exwaddcat=\exabcat^\op$
  and to a suitable \cooltech{}
  $\Dcool\excooltechdef$
  whose associated pointed category
  $\Dcool\localApof{\excooltech}=\connCh$
  is the shape of connective chain complexes;
  see \autoref{sec:example-Delta} for more details.
\item
  More generally,
  \autoref{thm:intro-main}
  specializes to the \inftycategorical{} \DKcorr{}
  originally sketched by Joyal~\cite[Section 35]{Joyal2008}
  \begin{equation}
    \Fun(\Dop,\exwaddcat)\xlra{\simeq}\connChof{\exwaddcat}
  \end{equation}
  between simplicial objects
  and coherent connective chain complexes
  in any weakly idempotent complete additive \inftycat{} $\exwaddcat$.
\item
  Denote by $\finsetp$ the category of finite pointed sets
  and by $\finsurj$ the category
  of (possibly empty) finite sets and surjections between them.
  Pirashvili\cite{Pirashvili2000} constructed an equivalence
  \begin{equation}
    \label{eq:intro-actual-pirashvili}
    \Fun(\finsetp,\Ab)\xlra{\simeq}\Fun(\finsurj,\Ab)
  \end{equation}
  between $\finsetp$-shaped
  and $\finsurj$-shaped diagrams\footnote{
    To be precise,
    Pirashvili only considers diagrams
    whose value on $\terminal\in\finsetp$
    and on $\emptyset\in\finsurj$ is zero;
    these diagrams correspond to each other
    under the equivalence~\eqref{eq:intro-actual-pirashvili}}
  of abelian groups.
  We recover this equivalence from
  \autoref{thm:intro-main}
  which more generally yields a natural equivalence
  \begin{equation}
    \label{eq:intro-pirashvili}
    \Fun(\finsetp,\exwaddcat)\xlra{\simeq}\Fun(\finsurj,\exwaddcat),
  \end{equation}
  between $\Gamma$-objects\footnote{
    Some authors define $\Gamma$ to be
    the category $\finsetp$ of finite pointed sets;
    we use Segal's original definition~\cite{Segal1974}
    which is dual, \ie $\Gamma\coloneqq\finsetpop$.
    Regardless of the convention,
    a $\Gamma$-object in $\exwaddcat$ is always a functor
    $\finsetp\to\exwaddcat$.
  } and $\finsurj$-shaped diagrams
  in any weakly idempotent complete preadditive\footnote{
    preadditive = has direct sums
  }
  \inftycat{} $\exwaddcat$;
  see \autoref{sec:cooltechs-from-factorizations} for more details.
\item
  Denote by $\FIsharp$
  the category of finite sets and partial injections;
  let $\finbij$
  be the groupoid of finite sets and bijections.
  For each commutative ground ring $\exring$,
  \cite[Theorem~4.1.5]{CEF2015}
  (which is a special case of \cite[Theorem~1.5]{Slominska2004})
  describes an equivalence
  \begin{equation}
    \label{eq:FIsharp-SG}
    \Fun(\FIsharp, \rMod{\exring})
    \xlra{\simeq}
    \Fun(\finbij,\rMod{\exring})
    \simeq
    \prod_{n\in\BN}\left(\lRep{\SG{n}}{\exring}\right)
  \end{equation}
  between the categories of $\FIsharp$-modules
  and of tuples of representations of all symmetric groups $\SG{n}$.
  Again, our main result
  generalizes this equivalence
  to coherent diagrams/representations with values in arbitrary
  weakly idempotent complete \pre{additive} \inftycats{}.
\item
  When $\exwaddcat$ is an
  idempotent complete additive \emph{ordinary} category,
  \autoref{thm:intro-main} recovers
  the general Dold--Kan type equivalence
  of \nameLS{}~\cite[Theorem~6.8]{LackStreet2015}
  which includes as special cases
  \eqref{eq:intro-DK-classical},
  \eqref{eq:intro-actual-pirashvili},
  \eqref{eq:FIsharp-SG}
  and many more.
  See \autoref{sec:comparison-Lack-Street} for a detailed comparison.

\item
  Some of the  equivalences
  of \autoref{thm:intro-main}---%
  including the one for $\Gamma=\finsetpop$ but not the one for $\Delta$---%
  were already established by Helmstutler~\cite{Helmstutler2014}
  in the language of model categories;
  see \autoref{rem:comparison-Helmstutler} for more details.
  Note that unlike \autoref{thm:intro-main},
  Helmstutler's result cannot be dualized so easily to yield,
  for instance,
  a model categorical version of the equivalence~\eqref{eq:intro-pirashvili}.
\item
  In a \emph{stable} \inftycat{} $\exstablecat$,
  coherent connective chain complexes can be encoded
  more conveniently as filtered objects,
  \ie as diagrams $\BN\to\exstablecat$;
  an explicit equivalence
  \begin{equation}
    \label{eq:Ariotta-identification}
    \Fun(\BN,\exstablecat)\simeq\connChof{\exstablecat}
  \end{equation}
  is part of \SAriotta{}'s \PhD thesis~\cite{Ariotta}.
  In this stable context,
  Lurie proved an
  \inftycategorical{} \DKcorr{}~\cite[Theorem 1.2.4.1]{Lurie2017}
  in the form of an equivalence
  \begin{equation}
    \label{eq:intro-Lurie-DK-functor}
    \Fun(\Dop,\exstablecat)\xra{\simeq}\Fun(\BN,\exstablecat),
  \end{equation}
  of \inftycats{};
  we expect this equivalence to agree with \eqref{eq:intro-DK-classical}
  under the identification \eqref{eq:Ariotta-identification}.
  Note that while both equivalences
  \eqref{eq:Ariotta-identification} and \eqref{eq:intro-Lurie-DK-functor}
  need the stability of $\exstablecat$
  to work,
  \autoref{thm:intro-main}---%
  just like the ordinary \DKcorr{}---%
  only needs that $\exwaddcat$ is \waddjectives{}.
  See \autoref{sec:comparison-Lurie} for a more detailled discussion.
\end{enumerate}

We now introduce the notion of a \cooltech{}
$\excooltechdef$
on which \autoref{thm:intro-main} is based.
It consists of a three-fold factorization system of type
\begin{equation}
  \bullet\xra{\epis}\bullet\xra{}\bullet\xra{\dualepis}\bullet,
\end{equation}
where the unnamed middle piece
together with suitably encoded zero relations
gives rise to the pointed category
$\localAp(\excooltech)$
appearing on the right side of the equivalence~\eqref{eq:intro-main}.
This notion is inspired by similar concepts appearing in
\cite{Slominska2004,Helmstutler2014,LackStreet2015}.
We give an illustration in the examples of $\Gamma$ and $\Delta$,
which are discussed in greater detail in
\autoref{sec:cooltechs-from-factorizations}
and \autoref{sec:example-Delta}.
\begin{itemize}
\item
  Every map $f\colon I\la J$ in $\Gamma=\finsetpop$
  can be written uniquely as the composition
  \begin{equation}
    \label{eq:intro-factorization-Gamma}
    I\lla\Im{f}\lla \frac{J}{\Ker f} \lla J,
  \end{equation}
  where
  \begin{itemize}
  \item
    the leftmost map is a bijection onto its image,
  \item
    the middle map is surjective
    and sends only the basepoint to the basepoint,
    in other words it just amounts to a surjection
    between the (possibly empty) sets
    obtained by omitting the basepoints,
  \item
    the rightmost map is bijective outside of its kernel
    (such maps are often called \buzzword{inert}).
  \end{itemize}
  The category of those arrows
  of which appear
  as the middle piece of \eqref{eq:intro-factorization-Gamma}
  is precisely (the opposite of)
  the category $\finsurj$;
  there are no zero relations in this case.
\item
  Every arrow in $\Delta$ can be written uniquely as the composition
  \begin{equation}
    \label{eq:intro-factorization-Delta}
    \bullet\xra{\codegmap{\geq 0}}
    \bullet\xra{(\cofacemap{0})}
    \bullet\xra{\cofacemap{>0}}
    \bullet
  \end{equation}
  where
  \begin{itemize}
  \item
    the left arrow $\codegmap{\geq 0}$
    is a (possibly empty) composition of codegeneracy maps,
  \item
    the middle arrow is either the identity or a $0$-th coface map,
  \item
    the right arrow is a (possibly empty) composition of
    $i$-th coface maps $\cofacemap{i}$ for $i>0$.
  \end{itemize}
  If one focuses only on the arrows of the second type,
  one obtains a category $\connCh$
  \begin{equation}
    \begin{tikzcd}[column sep=small]
      \numCh{0}
      \ar[rr,"\cofacemap{0}"]
      &
      \phantom{.}
      \ar[rr,bend right,dotted,no head]
      &\numCh{1}
      \ar[rr,"\cofacemap{0}"]
      &
      \phantom{.}
      \ar[rr,bend right,dotted,no head]
      &\numCh{2}
      \ar[rr,"\cofacemap{0}"]
      &
      \phantom{.}
      &\cdots
    \end{tikzcd}
  \end{equation}
  with zero relations
  by declaring the composite of two $0$-th coface maps to vanish
  (because it is not again a $0$-face map).
  Connective chain complexes are then exactly
  zero-preserving presheaves on $\connCh$.
  In order to properly encode the coherent zero relations
  in the \inftycategorical{} context,
  we actually consider $\connCh$ as a pointed category
  by adding an additional zero object through which all zero morphisms factor;
  for a more detailed explanation of this issue,
  see \autoref{sec:quotient-cats}.
\end{itemize}

\subsection{Acknowledgements}

This work was done during my Ph.D.~studies at the Hausdorff Center for
Mathematics (HCM);
I am very grateful to my supervisors Catharina Stroppel and Tobias Dyckerhoff
for supporting and encouraging me during this time.
Moreover, I would like to thank Gustavo Jasso
for many interesting discussions and for valuable feedback.
This research was funded by the Deutsche Forschungsgemeinschaft
(DFG, German Research Foundation)
under Germany's Excellence Strategy - GZ 2047/1, Projekt-ID 390685813.
This paper was written while
\theauthor{} was a guest at the Max Planck Institute for Mathematics (MPIM).

\section{Preliminaries}

\subsection{\inftycategorical{} notation and tools}

\label{sec:infty-cat-tools}
\newcommand\excat[1][]{\manyprime[#1]{C}}
\newcommand{\exdiagram}[1][]{\manyprime[#1]{D}}
\newcommand{\locald}{d}
\newcommand\localc[1][]{\manyprime[#1]{c}}

Throughout this article we use the language of \inftycats{}/quasi-categories
as developed by Joyal and Lurie;
our main references are \cite{Lurie2009} and \cite{Cisinski2019}.

We view the theory of \inftycats{} as an extension of ordinary category theory
by identifying an ordinary category $\excat$ with its nerve $\nerve{\excat}$.
Typically, we use ordinary capital letters
(\eg $\excat,\exdiagram,\exptdordcat,\exabcat$)
for $1$-categories
and the corresponding Euler Script letters
(\eg $\exinftycat,\inftycatcal{\exdiagram}, \exptdcat,\exwaddcat$)
for \inftycats{}.

We write $\Fun(\exinftycat[1],\exinftycat)$
for the \inftycat{} of functors $\exinftycat[1]\to\exinftycat$.
Given a small category $\exdiagram$ and an \inftycat{} $\extargetcat$
a diagram of shape $\exdiagram$ with values in $\extargetcat$
is a functor $\exdiagram\to\extargetcat$;
a $\extargetcat$-valued presheaf on $\exdiagram$
is a functor $\exdiagram^\op\to\extargetcat$.

Given a category $\excat$
and an object $\localc\in\excat$,
we denote by
$\overcat{\excat}{\localc}$
the category of objects over $\localc$;
its objects are arrows of the form $\bullet\to\localc$.
Given a functor $\excat[1]\to\excat$ of categories
and an object $\localc\in\excat$,
we denote by
$\overcat[\excat]{\excat[1]}{\localc}$
the relative overcategory
defined as the fiber product
$\excat[1]\fiberproduct{\excat}\overcat{\excat}{\localc}$
(with the functor $\excat[1]\to\excat$ left implicit).
Dually, the symbols
$\undercat{\localc}{\excat}$
and
$\undercat[\excat]{\localc}{\excat[1]}$
denote absolute and relative undercategories.

We denote by
$\rcone{\exinftycat}\coloneqq\exinftycat\join\set{\conemax}$
and $\lcone{\exinftycat}\coloneqq\set{\conemin}\join\exinftycat$
the \inftycats{} obtained from the \inftycat{} $\exinftycat$
by adjoining a new terminal object $\conemax$
or initial object $\conemin$, respectively.
We denote by $\pi_0\exinftycat$
the set of equivalence classes of objects
(\ie the $0$-truncation),
by $\htpycat{\exinftycat}$ the homotopy category
(\ie the $1$-truncation)
and by $\grpdcore{\exinftycat}$
the groupoid core
(obtained by discarding \non{invertible arrows})
of $\exinftycat$.

The main tool of this article is the theory
of \inftycategorical{} (co)limits and Kan extensions
as developed in \cite[Chapter 4]{Lurie2009}.
We recall briefly the following key statements,
which are analogs from classical facts of ordinary category theory
and will be used throughout this article
without further mention:

\begin{itemize}
\item\cite[Definition~4.3.2.2]{Lurie2009}
  Right/left Kan extension along a fully faithful functor
  $\excat[1]\hra\excat$
  can be computed and characterized pointwise at each
  $\locald\in\excat$
  by the induced limit/colimit of shape
  $\undercat[\excat]{d}{\excat[1]}$
  and
  $\overcat[\excat]{\excat[1]}{d}$,
  respectively.
\item\cite[Proposition~4.3.2.15]{Lurie2009}
  Restriction along a fully faithful functor
  $\excat[1]\hra\excat$
  induces an equivalence of \inftycats{}
  between the full \sub{categories} of
  $\Fun(\excat,\extargetcat)$
  and
  $\Fun(\excat[1],\extargetcat)$
  consisting of those functors which \emph{are} a right/left Kan extension
  and those functors which \emph{have} a right/left Kan extension,
  respectively.
\item\cite[Corollary 4.3.2.16, Proposition 4.3.2.17]{Lurie2009}
  If every functor $\excat[1]\to\extargetcat$
  admits a right/left Kan extension along
  the fully faithful functor $\excat[1]\hra\excat$
  then there is a unique fully faithful right/left Kan extension functor
  $\Fun(\excat[1],\extargetcat)
  \to
  \Fun(\excat,\extargetcat)$
  which is right/left adjoint to the restriction functor;
  its essential image is spanned by those functors
  $\excat\to\extargetcat$
  which are a right/left Kan extension along
  $\excat[1]\hra\excat$.
\item\cite[Proposition~4.1.3.1]{Lurie2009}
  A functor
  $\exdiagram[1]\to\exdiagram$
  between ordinary categories is homotopy terminal\footnote{
    Joyal and Lurie would say \buzzword{cofinal}
    which, confusingly, is the word Cisinski uses for the dual concept
    (what we call homotopy initial).
    We avoid this potential confusion by using
    the hopefully unambiguous terminology of Dugger~\cite{Dugger}.
  }
  if and only if
  each undercategory
  $\undercat[\exdiagram]{d}{\exdiagram[1]}$
  (for each $\locald\in\exdiagram$)
  is \wcontractible{},
  \ie has contractible geometric realization
  $\geomrealization{
    \nerve{\undercat[\exdiagram]{d}{\exdiagram[1]}}
  }$.
  Dually $\exdiagram[1]\to\exdiagram$
  is homotopy initial
  if and only if
  each overcategory
  $\overcat[\exdiagram]{\exdiagram[1]}{d}$
  is \wcontractible{}.
\item\cite[Proposition~4.1.1.8]{Lurie2009}
  The limit/colimit of a $\exdiagram$-shaped diagram
  $\exdiagram\to\extargetcat$
  can be computed after \pre{composing}
  with any homotopy initial/terminal functor
  $\exdiagram[1]\to\exdiagram$.
\end{itemize}

\newcommand{\localW}{W}

Finally, we remind the reader that a localization\footnote{
  Here our terminology differs from Lurie's
  who uses the word \qquote{localization}
  to refer to a special kind of localization functor
  which admits a fully faithful right adjoint
  (see \cite[Definition~5.2.7.2 and Warning~5.2.7.3]{Lurie2009}).
} of an \inftycat{}
$\exinftycat$ at a class $\localW\subset\exinftycat$ of arrows
is a functor
$\exinftycat\to\localizedat{{\exinftycat}}{\localW}$
which is universal amongst all functors that
send the arrows in $\localW$ to equivalences.
More precisely,
for each \inftycat{} $\exinftycat[1]$
the restriction functor
\begin{equation}
  \Fun(\localizedat{{\exinftycat}}{\localW},\exinftycat[1])
  \lra
  \Fun(\exinftycat,\exinftycat[1])
\end{equation}
is fully faithful with essential image consisting of those functors
$\exinftycat\to\exinftycat[1]$
that send all arrows in $\localW$ to equivalences in $\exinftycat[1]$
(see, for example \cite[Definition~7.1.2
]{Cisinski2019}.
Such \inftycategorical{} localizations always exist
and are essentially unique, see \cite[Proposition~7.1.3]{Cisinski2019}.

\subsection{Pointed \inftycats{}}

\newcommand{\localx}{x}
\newcommand{\localy}{y}
\newcommand{\localf}{f}
\newcommand{\localX}{X}
\newcommand{\localY}{Y}
\newcommand{\localK}{K}
\newcommand{\exconediagram}{D}

Recall, that an \inftycat{} $\exptdcat$ is called \introduce{pointed}
if it has a zero object,
\ie an object $\zerobj\in\exptdcat$
which is both initial and terminal in $\exptdcat$.
A functor $\exptdcat[1]\to\exptdcat$ between pointed \inftycats{}
is called \introduce{pointed}
if it sends one (equivalently, each) zero object of $\exptdcat[1]$
to a zero object of $\exptdcat$.
We denote by $\Catinftyptd$ the \inftycat{}
of (small) pointed \inftycats{} and pointed functors between them;
it comes equipped with a canonical forgetful functor
\begin{equation}
  \Catinftyptd \lra \Catinfty.
\end{equation}
Given two pointed \inftycats{} $\exptdcat[1]$ and $\exptdcat$,
we denote by
$\ptdFun{\exptdcat[1]}{\exptdcat}\subset \Fun(\exptdcat[1],\exptdcat)$
the full \sub{category} spanned by the pointed functors.

\begin{Cstr}[Free pointed category]\label{cstr:ordinary-free-base-pt}
  Let $\exordcat$ be an ordinary category.
  We define a pointed category $\addbasept\exordcat$
  by freely adjoining a zero object to $\exordcat$.
  Explicitly, it is described as follows:
  \begin{itemize}
  \item
    The objects of $\addbasept\exordcat$ are the objects of $\exordcat$
    plus an additional object $\zerobj$.
  \item
    For every object $\localx\in\addbasept\exordcat$
    we put
    \begin{equation}
      \addbasept\exordcat(\localx,\zerobj)=\set{\zeromap}
      \intxt{and}
      \addbasept\exordcat(\zerobj,\localx)=\set{\zeromap}
    \end{equation}
    (in other words, $\zerobj\in\addbasept\exordcat$
    is a zero object as the notation suggests).
    Given objects $\localx,\localy \in\exordcat$,
    we set
    \begin{equation}
      \addbasept\exordcat(\localx,\localy) \coloneqq
      \exordcat(\localx,\localy)\disjunion\set{\zeromap}
    \end{equation}
    where here $\zeromap$ denotes the composite map
    $\localx\to\zerobj\to \localy$.
  \item
    The composition in $\addbasept\exordcat$
    is inherited from the composition in $\exordcat$.
    \qedhere
  \end{itemize}
  The pointed category $\addbasept\exordcat$ comes equipped
  with the canonical (\non{full}) inclusion functor
  $\exordcat\to \addbasept\exordcat$.
\end{Cstr}

\begin{Cstr}[Free pointed \inftycat{}]
  \label{cstr:infty-free-base-pt}
  Let $\exinftycat$ be an \inftycat{}.
  Denote by
  \begin{equation}
    \lrcone{\exinftycat}
    \coloneqq
    \set{\conemin}\join\exinftycat \join \set{\conemax}
  \end{equation}
  the \inftycat{} obtained from $\exinftycat$
  by freely adjoining an initial object $\conemin$
  and a terminal object $\conemax$.
  We define $\addbasept{\exinftycat}$
  to be the localization of $\lrcone\exinftycat$
  at the (essentially unique) edge $\conemin\to\conemax$
  connecting the initial to the terminal object.
  The \inftycat{} $\addbasept{\exinftycat}$ is pointed
  (since localizations preserve both initial and terminal objects\footnote{
    This follows, for instance, from Proposition 7.1.10 in \cite{Cisinski2019}})
  and comes equipped with the defining functor
  $
  \exinftycat
  \hra
  \lrcone{\exinftycat}
  \to
  \addbasept{\exinftycat}.
  $
\end{Cstr}

If the category $\exinftycat$ in \autoref{cstr:infty-free-base-pt}
happens to be an ordinary category,
then $\lrcone{\exinftycat}$ is again an ordinary category.
It is however not clear \latin{a priori} that
the the same is true for $\addbasept\exinftycat$,
because the localization procedure has the potential to turn an ordinary
category into one that isn't.
The following lemma addresses this issue.

\begin{Lem}
  \label{lem:compare-add-basept-ord}
  Let $\exordcat$ be an ordinary category.
  Then the functor $\exordcat\to \addbasept\exordcat$
  from \autoref{cstr:infty-free-base-pt} agrees
  with the one from \autoref{cstr:ordinary-free-base-pt}.
  In particular, $\addbasept\exordcat$ is an ordinary category again.
\end{Lem}

\begin{Prf}
  \newcommand{\localga}{\gamma}
  \newcommand{\localmm}{\spanby{\conemin,\conemax}}
  \newcommand{\localdimxpart}[1]{n(#1)}
  \newcommand{\localsimplexC}{\localx}
  \newcommand{\localintersNA}{t}
  \newcommand{\localinters}[1]{\localintersNA(#1)}
  \newcommand{\localk}{k}
  \newcommand{\locali}{i}
  \newcommand{\localbigdim}{m}
  \newcommand{\localzerotpl}[1]{\zerobj^{#1}}
  \newcommand{\locald}{d}
  \newcommand{\localfiltrCp}[1]{\nerve{\addbasept{\exordcat}}^{\leq #1}}
  \newcommand{\localbigsimplex}{\sigma(\localk,\localsimplexC,\localintersNA)}
  Let $\addbasept\exordcat$ be as in \autoref{cstr:ordinary-free-base-pt}
  and consider the canonical functor
  \begin{equation}
    \localga\colon
    \lrcone{\exordcat}
    =
    \set{\conemin}\join\exordcat\join\set{\conemax}
    \lra\addbasept\exordcat
  \end{equation}
  given by the canonical inclusion of $\exordcat$
  and by $\conemin,\conemax\mapsto\zerobj$.
  We need to show that $\localga$ exhibits $\addbasept\exordcat$
  as the \inftycategorical{} localization of
  $\set{\conemin}\join\exordcat \join \set{\conemax}$
  at the unique map $\conemin\to\conemax$.
  Denote by $\localmm$ the full \sub{category} of
  $\lrcone{\exordcat}$
  spanned by $\conemin$ and $\conemax$.
  Since $\localmm\cong\simplex 1$ is weakly contractible,
  it follows by comparing universal properties
  that the desired localization can be computed as the pushout
  $\lrcone{\exordcat} \pushout{\localmm}{\set{\zerobj}}$
  (of \inftycats{}).
  To conclude the proof,
  it therefore suffices to show that---after passing to nerves---%
  the canonical square
  \begin{equation}
    \cdsquareOpt
    {\localmm}
    {\set{\zerobj}}
    {\lrcone{\exordcat}}
    {\addbasept\exordcat}
    {}{hookrightarrow}{}{}
  \end{equation}
  of categories becomes a (Joyal) homotopy pushout of simplicial sets.
  Since the left vertical map is a monomorphism,
  it suffices to show that the map
  \begin{equation}
    \label{eq:inprf:map-pushout-to-Cp}
    {\nerve{\set{\zerobj}}}
    \pushout {\nerve{\localmm}}
    {\nerve{\lrcone{\exordcat}}}
    \lra\nerve{\addbasept\exordcat}
  \end{equation}
  from the (strict) pushout of simplicial sets is a (Joyal) weak equivalence;
  we will now show that it is in fact
  an inner anodyne extension.

  The simplices of
  $\addbasept{\exordcat}$
  can be described explicitly as follows:
  Each $\localbigdim$-simplex of $\nerve{\addbasept{\exordcat}}$
  is of the form
  \begin{equation}
    \label{eq:inprf:bigchain-many-zeros}
    \localbigsimplex\colon
    \localzerotpl{\localinters{0}}
    \to
    \localsimplexC^1
    \to
    \localzerotpl{\localinters{1}}
    \to
    \localsimplexC^2
    \to
    \localzerotpl{\localinters{2}}
    \to
    \cdots
    \to
    \localzerotpl{\localinters{\localk-1}}
    \to
    \localsimplexC^{\localk}
    \to
    \localzerotpl{\localinters{\localk}},
  \end{equation}
  where
  \begin{itemize}
  \item
    $\localk$ is a natural number
  \item
    each
    $
    \localsimplexC^\locali\colon
    \localsimplexC^{\locali}_{0}
    \to\dots\to
    \localsimplexC^{\locali}_{\localdimxpart{\locali}}
    $
    (for $1\leq\locali\leq\localk$)
    is an $\localdimxpart{\locali}$-simplex
    of $\nerve\exordcat$.
  \item
    $
    \localinters{0},
    \dots,
    \localinters{\localk}
    $
    are natural numbers
    of which all but $\localinters{0}$ and $\localinters{\localk}$
    are required to be positive.
  \item
    $\localzerotpl{\localinters{\locali}}$
    denotes a chain
    $\zerobj\to\cdots\to\zerobj$
    with $\localinters{\locali}$ many zeros.
  \item
    the dimension
    $
    \localbigdim
    \coloneqq
    \localinters{0}
    -1
    +
    \sum_{\locali=1}^{\localk}
    \left(
      \localdimxpart{\locali}+1+\localinters{\locali}
    \right)
    $
    is \non{negative}.
  \end{itemize}
  Denote by
  $\localfiltrCp{\locald}\subset\nerve{\addbasept{\exordcat}}$
  the simplicial \sub{set} containing those simplices
  $\localbigsimplex$ with $\localk\leq\locald$.
  The following are straightforward to verify:
  \begin{enumerate}
  \item
    \label{item:wanted-pushout-is-Cp1}
    The map~\eqref{eq:inprf:map-pushout-to-Cp}
    induces an isomorphism
    $
    {\nerve{\set{\zerobj}}}
    \pushout {\nerve{\localmm}}
    {\nerve{\lrcone{\exordcat}}}
    \xra{\cong}
    \localfiltrCp{1}
    $.
  \item
    \label{item:pushout-Cpd-Cpd-1}
    For each $\locald\geq 1$,
    we have a pushout of simplicial sets
    \begin{equation}
      \label{eq:inpf-pushout-for-Cpd}
      \cdsquareOpt[cC]
      {
        \coprod\limits_{\localk,\localsimplexC,\localintersNA}
        {
          \simplex{
            \localinters{0}
            +
            \localdimxpart{1}
            +
            \localinters{1}
          }
          \pushout{\simplex{\localinters{1}-1}}
          \simplex{\localbigdim'}
        }
      }
      {\localfiltrCp{\locald-1}}
      {
        \coprod\limits_{\localk,\localsimplexC,\localintersNA}
        \simplex{\localbigdim}
      }
      {\localfiltrCp{\locald}}
      {}
      {hookrightarrow}
      {hookrightarrow}
      {"(\localbigsimplex)"}
    \end{equation}
    which corresponds to the decomposition of
    each chain~\eqref{eq:inprf:bigchain-many-zeros}
    into the two overlapping chains
    \begin{equation}
      \localzerotpl{\localinters{0}}
      \to
      \localsimplexC^1
      \to
      \localzerotpl{\localinters{1}}
      \intxt{and}
      \localzerotpl{\localinters{1}}
      \to
      \localsimplexC^2
      \to
      \localzerotpl{\localinters{2}}
      \to
      \cdots
      \to
      \localzerotpl{\localinters{\localk-1}}
      \to
      \localsimplexC^{\localk}
      \to
      \localzerotpl{\localinters{\localk}}
    \end{equation}
    of dimensions
    $
    \localinters{0}
    +
    \localdimxpart{1}
    +
    \localinters{1}
    $
    and
    $\localbigdim'\coloneqq
    -1
    +\localinters{1}
    +
    \sum_{\locali=2}^{\localk}
    \left(
      \localdimxpart{\locali}+1+\localinters{\locali}
    \right),
    $
    respectively.
  \item\label{item:Cp-isfilteredunion}
    The simplicial set
    $\nerve{\addbasept{\exordcat}}$
    is the union of the ascending chain
    $
    \localfiltrCp{1}
    \subset
    \localfiltrCp{2}
    \subset
    \cdots
    $
    of simplicial \sub{sets}.
  \end{enumerate}
  The left vertical map in the square \eqref{eq:inpf-pushout-for-Cpd}
  is an inner anodyne extension;
  it follows from
  \ref{item:wanted-pushout-is-Cp1},
  \ref{item:pushout-Cpd-Cpd-1} and
  \ref{item:Cp-isfilteredunion}
  that the same is true for the map
  \eqref{eq:inprf:map-pushout-to-Cp};
  this concludes the proof.
\end{Prf}

\begin{Rem}
  In view of \autoref{lem:compare-add-basept-ord},
  we are justified in tacitly assuming that the free pointed \inftycat{}
  $\addbasept\exordcat$ on an ordinary category $\exordcat$
  is given by the explicit description of \autoref{cstr:ordinary-free-base-pt}.
\end{Rem}

The following lemma establishes the universal property of
the free pointed \inftycat{} construction.

\begin{Prop}
  \label{prop:univ-prop-add-basept}
  Let $\exinftycat$ be a (small) \inftycat{}.
  For every pointed \inftycat{} $\exptdcat$,
  restriction along the functor $\exinftycat\to\addbasept{\exinftycat}$
  induces an equivalence
  \begin{equation}
    \ptdFun{\addbasept{\exinftycat}}{\exptdcat}
    \xra{\simeq}
    \Fun(\exinftycat,\exptdcat).
  \end{equation}
  of \inftycats{}.
  In particular, the construction
  $\exinftycat\mapsto\addbasept\exinftycat$
  yields a left adjoint to the forgetful functor
  $\Catinftyptd\to\Catinfty$
\end{Prop}

\begin{Prf}
  \def\tempa{1}
  \def\tempb{2}
  \def\tempc{3}
  \def\tempd{4}
  \def\tempe{5}
  \def\tempf{6}
  The functors
  $
  \exinftycat
  \hra
  \lcone\exinftycat
  \hra
  \lrcone\exinftycat
  \lra
  \addbasept\exinftycat
  $
  induce the following commutative diagram of functor \inftycats{}
  and their various \sub{categories} defined as indicated:
  \begin{equation}
    \begin{tikzcd}
      \Fun(\addbasept\exinftycat,\exptdcat)
      \ar[ddd]\ar[r,hookleftarrow]
      \ar[rdd,"{\simeq}","{\tempa}"']
      &
      \ptdFun{\addbasept\exinftycat}{\exptdcat}
      \ar[rdd,"{\simeq}","{\tempb}"']
      \\
      \\
      &
      \set{\text{inverts }\conemin\to\conemax}
      \ar[r,hookleftarrow]
      &
      \set{
        \text{inverts }
        \conemin\to\conemax
        \text{ and }
        \conemin,\conemax\mapsto\zerobj
      }
      \ar[from=dd,"=","{\tempf}"',hookrightarrow]
      \\
      \Fun(\lrcone{\exinftycat},\exptdcat)
      \ar[ddd]\ar[rd,hookleftarrow]\ar[ru,hookleftarrow]
      \\
      &
      \set{\conemax\mapsto\zerobj}
      \ar[r,hookleftarrow]
      &
      \set{\conemin,\conemax\mapsto\zerobj}
      \\
      \\
      \Fun(\lcone{\exinftycat},\exptdcat)
      \ar[from=ruu,"{\simeq}"',"{\tempc}"]
      \ar[d]\ar[r,hookleftarrow]
      &
      \set{\conemin\mapsto \zerobj}
      \ar[from=ruu,"{\simeq}"',"{\tempd}"]
      \\
      \Fun(\exinftycat,\exptdcat)
      \ar[from=ru,"{\simeq}"',"{\tempe}"]
    \end{tikzcd}
  \end{equation}
  Restriction along
  $\lrcone\exinftycat\to\addbasept\exinftycat$
  induces the equivalence \tempa{}
  by the universal property of the localization.
  The functors labeled by \tempc{} and \tempe{} are equivalences
  because they have an inverse given
  by right Kan extension and left Kan extension, respectively
  (using that $\zerobj\in\exptdcat$
  is a terminal and an initial object, respectively).
  The equivalences \tempb{} and \tempd{} are induced
  by restricting to appropriate full \sub{categories}.
  Since $\Map_{\exptdcat}(\zerobj,\zerobj)\simeq\pt$,
  every functor $\lrcone\exinftycat\to\exptdcat$ which sends
  $\conemin$ and $\conemax$ to zero objects
  must invert the edge $\conemin\to\conemax$;
  thus the inclusion labeled \tempf{}
  is an equality of full \sub{categories}.
  The result follows.
\end{Prf}

\begin{Lem}
  \label{lem:factors-of-equis-are-equis}
  Let $\exptdcat$ be a pointed \inftycat{} and let
  $\setP{f_i\colon \localX_i\to \localY_i}{i\in I}$
  be a finite set of morphisms in $\exptdcat$.
  Assume that the product
  \begin{equation}
    \prod_{i\in I}f_i\colon
    \prod_{i\in I}\localX_i\lra
    \prod_{i\in I}\localY_i
  \end{equation}
  exists in $\exptdcat$ and is an equivalence.
  Then for each $i\in I$,
  the morphism $f_i\colon \localX_i\to \localY_i$ is an equivalence.
\end{Lem}

\begin{Prf}
  Given an inverse $g\colon\prod_i \localY_i\to \prod_i\localX_i$ to $\prod f_i$,
  it is easy to see that for each $j\in I$ the composition
  \begin{equation}
    \label{eq:diagonal-matrix-inverse}
    \localY_j\xra{\factorincl{j}}\prod_i \localY_i\xra{g}\prod_i \localX_i\xra{\factorproj{j}}\localX_j,
    \intxt{where}
    \factorproj{i}\factorincl{j}\coloneqq
    \begin{cases}
      \Id\colon \localY_j\to \localY_j & \text{if } i=j\\
      \zeromap\colon \localY_j\to \localY_i & \text{if } i\neq j
    \end{cases}
  \end{equation}
  is an inverse of $f_j$.
\end{Prf}

In a pointed \inftycat{} it makes sense to talk about fibers and cofibers
which are the \inftycategorical{} analog of kernels and cokernels.
The \introduce{fiber} and \introduce{cofiber}
of an arrow $f\colon\localX\to\localY$
are the pullback and pushout
of the diagrams
\begin{equation}
  \begin{tikzcd}
    &\zerobj\ar[d]\\
    \localX\ar[r,"f"]&\localY
  \end{tikzcd}
  \intxt{and}
  \begin{tikzcd}
    \localX\ar[d]\ar[r,"f"]&\localY\\
    \zerobj
  \end{tikzcd}
\end{equation}
respectively.
More generally,
we define the \introduce{total cofiber}
$\totfib{\exconediagram}$
of a conical diagram
$\exconediagram\colon\rcone{\localK}\to \exptdcat$
as the cofiber of the canonical map
$
\colim (\localK\subset\rcone{\localK}\xra{\exconediagram}\exptdcat)
\to \exconediagram(\conemax)
$
and the \introduce{total fiber}
$\totcof{\exconediagram}$
of a conical diagram
$\exconediagram\colon\lcone{\localK}\to \exptdcat$
as the fiber of the canonical map
$
\exconediagram(\conemin)
\to
\lim (\localK\subset\lcone{\localK}\xra{\exconediagram}\exptdcat)
$.
To recover the case of the ordinary fiber/cofiber
set $\localK=\simplex 0$,
hence $\rcone{\localK}\cong\simplex 1\cong\lcone{\localK}$.

Another way of computing the total cofiber
(\resp total fiber)
of a $\rcone{\localK}$-shaped
(\resp $\lcone{\localK}$-shaped)
diagram $\exconediagram$
is to first pass to its right (\resp left) Kan extension
along the first inclusion
$\rcone{\localK}\hra\rcone{\localK}\pushout{\localK}\rcone{\localK}$
(\resp
$\lcone{\localK}\hra\lcone{\localK}\pushout{\localK}\lcone{\localK}$)---%
which is given explicitly by setting
the value on the cone point of the second copy of
$\rcone{\localK}$ (\resp $\lcone{\localK}$)
to $\zerobj\in\exptdcat$---%
and then taking the colimit (\resp limit) of this diagram.
The advantage of this description is that it is well defined
even if the colimit (\resp limit)
of $\restr\exconediagram\localK$ does not exist in $\exptdcat$.

\subsection{Quotient categories and coherent chain complexes}

\label{sec:quotient-cats}
\newcommand{\exideal}{S}
\newcommand{\localx}{x}
\newcommand{\localy}{y}
\newcommand{\localf}{f}
\newcommand{\exquotient}{{\quotcat{\exordcat}{\exideal}}}

A chain complex in an ordinary pointed category $\exptdordcat$
is a diagram $\BZ^\op\to\exptdordcat$,
which we might depict as
\begin{equation}
  \cdots\xla{\exdiff}
  \bullet
  \xla{\exdiff}
  \bullet
  \xla{\exdiff}
  \bullet
  \xla{\exdiff}\cdots
\end{equation}
such that any composite of more than one $\exdiff$
is sent to the zero morphism in $\exptdordcat$.
In other words,
the category of chain complexes in $\exptdordcat$
is a full subcategory of the category
$\Fun(\BZ^\op,\exptdordcat)$
of $\exptdordcat$-valued \pre{sheaves} on $\BZ$.
In the \inftycategorical{} world,
this naive definition would no longer be satisfactory because
\begin{itemize}
\item
  for a map in an \inftycat{}, being zero is no longer a \emph{property}
  but the \emph{structure} of an explicit null-homotopy and
\item
  there should be be higher coherence data
  exhibiting all the trivializations
  $\exdiff\circ\dots\circ\exdiff\simeq\zeromap$
  as compatible
\end{itemize}
Let $\exordcat$ be a category
equipped with an ideal $\exideal\subseteq\exordcat$
(\ie a set of arrows satisfying
$\exordcat\circ\exideal\circ\exordcat\subseteq\exideal$),
we would like to say what it means to equip a diagram
$\exordcat\to\exptdcat$
with a coherent trivialization of all arrows in $\exideal$.

\begin{Cstr}
  \label{cstr:quotient-pointed-cat}
  We define a pointed category $\exquotient$ as follows:
  \begin{itemize}
  \item
    The objects of $\exquotient$
    are the objects $\localx\in\exordcat$
    plus an additional zero object $\zerobj$.
  \item
    The morphisms of $\exquotient$ are determined by setting
    \begin{equation}
      \exquotient(\localx,\localy)
      \coloneqq
      \frac{\exordcat(\localx,\localy)}{\exideal}
      \cong
      \setP{\localf\in\exordcat(\localx,\localy)}{\localf\notin\exideal}
      \disjunion\set{\localx\to\zerobj\to\localy}
    \end{equation}
    for $\localx,\localy\in\exordcat$,
    with composition induced by the one in $\exordcat$.
  \end{itemize}
  The category $\exquotient$
  comes equipped with the canonical functor
  $\exordcat\to\exquotient$
  which is the identity on objects
  and sends precisely the arrows in $\exideal$ to zero.
\end{Cstr}

\begin{Rem}
  If $\localx\in\exordcat$ is an object with $\Id_\localx\in\exideal$
  then the unique morphisms $\localx\to\zerobj$ and $\zerobj\to\localx$
  are mutually inverse isomorphisms in $\exquotient$.
\end{Rem}

\begin{Def}
  Let $\exordcat\to\exptdcat$ be a $\exordcat$-shaped diagram
  in a pointed \inftycat{} $\exptdcat$.
  We say that a \introduce{trivialization} of all arrows in $\exideal$
  is an extension of
  $\exordcat\to\exptdcat$
  along $\exordcat\to\exquotient$
  to a pointed functor
  $\exquotient\to\exptdcat$.
\end{Def}

\begin{Expl}
  \begin{itemize}
  \item
    The quotient $\quotcat{\exordcat}{\emptyset}$
    of $\exordcat$ by the empty ideal
    is the free pointed category $\addbasept{\exordcat}$.
    Hence \autoref{prop:univ-prop-add-basept} can be read as saying that
    the empty set of arrows can always be trivialized in a unique way.
  \item
    If the category $\exordcat$ is already pointed
    and $\exideal=\zeroideal$ consists of all zero maps
    $\bullet\to\zerobj\to\bullet$
    then $\quotcat{\exordcat}{\zeroideal}\cong\exordcat$.
  \item
    Every category $\exordcat$ has an ideal
    consisting of all \non{isomorphisms};
    the corresponding quotient
    $\quotcat{\exordcat}{\not\simeq}$
    is the free pointed category $\quotbyneq{\exordcat}$
    on the groupoid core $\exordcat^\simeq$ of $\exordcat$.
    \qedhere
  \end{itemize}
\end{Expl}

\begin{Def}
  We denote by $\Ch\coloneqq\quotcat{\BZ}{\idealtwodiff}$
  the quotient of the poset $\BZ$
  by the ideal $\idealtwodiff$ of all maps
  $\numCh{\nn}\ra\numCh{\mm}$
  with $\mm-\nn\geq 2$.
  A \introduce{coherent chain complex} in $\exptdcat$
  is a pointed pre{sheaf} $\Ch^\op\to\exptdcat$;
  we denote by
  $\Chof{\exptdcat}\coloneqq\ptdFun{\Ch^\op}{\exptdcat}$
  the \inftycat{} of coherent chain complexes in $\exptdcat$.
  Similarly, we set
  $\connCh\coloneqq\quotcat{\BN}{\idealtwodiff}$
  and define the \inftycats{} of
  \introduce{connective chain complexes}
  in $\exptdcat$ as
  $\connChof{\exptdcat}\coloneqq\ptdFun{\connCh^\op}{\exptdcat}$.
\end{Def}

\begin{Rem}
  \label{rem:1-trivializing-ff}
  If $\exptdordcat$ is a pointed $1$-category
  then it is straightforward to check that the restriction functor
  \begin{equation}
    \ptdFun{\exquotient}{\exptdordcat}
    \lhra
    \Fun(\exordcat,\exptdordcat)
  \end{equation}
  is fully faithful and that
  the essential image consists of those functors
  $\exordcat\to\exptdordcat$
  which send arrows in $\exideal$ to zero maps in $\exptdordcat$.
  This means that
  \qquote{sending arrows in $\exideal$ to zero}
  is a \emph{property}
  which a diagram $\exordcat\to\exptdordcat$ might or might not have.
  If $\exptdcat$ is an \inftycat{}, this is no longer true in general:
  specifying a lift of a diagram
  $\exordcat\to\exptdcat$
  to a pointed diagram
  $\exquotient\to\exptdcat$
  might require an infinite amount of additional \emph{structure}.
\end{Rem}

\begin{Rem}
  \label{rem:what-about-pointed-enriched}
  Another way to make the notion
  of trivialization of arrows in $\exideal$ precise
  would have been to work with \inftycats{} enriched in \emph{pointed}
  spaces or even in pairs of spaces.
  Then we could study pairs-enriched diagrams
  $\exordcat\to\exptdcat$,
  where $\exordcat$ is pairs-enriched via $\exideal$
  and where $\exptdcat$ is pairs-enriched
  (even $\SpacesP$-enriched)
  via the zero maps.
  From this perspective one can see in a different way
  how the  additional structure encoded in such trivializations comes in:
  unlike the forgetful functor
  $\SetP\to\Set$
  from pointed sets to sets,
  the \roughly{forgetful} functor
  $\SpacesP\to\Spaces$
  from the \inftycat{} of pointed spaces to the \inftycat{} of spaces
  is not faithful
  and in fact not even injective on $\pi_0$ of mapping spaces.
\end{Rem}

\subsection{Additive and \pre{additive} \inftycats{}}

\label{sec:waddcats}
\newcommand{\localx}[1][]{\manyprime[#1]{X}}
\newcommand{\localy}[1][]{\manyprime[#1]{Y}}

\begin{Def}\cite[Definitions~2.1 and 2.6]{GGN2015}
  \label{def:additive-cat}
  An \inftycat{} $\exwaddcat$ with finite products and coproducts
  is called \introduce{\pre{additive}} if
  \begin{itemize}
  \item
    it is pointed, \ie the canonical map
    $\initial\xra{\simeq}\terminal$
    from the initial objects to the terminal object is an equivalence.
  \item
    for any two objects $\localx,\localx[1]\in\exwaddcat$,
    the canonical morphism
    \begin{equation}
      \label{eq:coproduct-to-product}
      \left(
        \scriptsize{
          \begin{matrix}
            \Id&\zeromap\\
            \zeromap&\Id
          \end{matrix}
        }
      \right)
      \colon
      \localx\sqcup \localx[1]
      \xra{\simeq}
      \localx\times \localx[1]
    \end{equation}
    (which exists, since $\exwaddcat$ is pointed)
    is an equivalence.
  \end{itemize}
  The \inftycat{} $\exwaddcat$ is called \introduce{additive} if additionally
  \begin{itemize}
  \item
    for each object $\localx\in\exwaddcat$,
    the \buzzword{shear map}
    \begin{equation}
      \label{eq:shear-map}
      \left(
       \scriptsize{
          \begin{matrix}
            \Id&\Id\\
            \zeromap&\Id
          \end{matrix}
        }
      \right)
      \colon
      \localx\sqcup \localx
      \xra{\simeq}
      \localx\times \localx
    \end{equation}
    is an equivalence.
  \end{itemize}
  A functor between \pre{additive} \inftycats{} is called
  \introduce{additive}
  if it is pointed and preserves finite products
  (or, equivalently, finite coproducts).
\end{Def}

\begin{Rem}
  Since products and coproducts in a \pre{additive} \inftycat{}
  are canonically identified,
  it is customary to call them \introduce{direct sums}
  and denote them by the same symbol $\oplus$.
\end{Rem}

\begin{Rem}
  When specializing to the case where $\exwaddcat$ is an ordinary category,
  \autoref{def:additive-cat} recovers the classical notion
  of an additive category
  (as defined, for instance, in \cite[Chapter~VIII]{MacLane1998}).
  However, we warn the reader that our use of the word \qquote{\pre{additive}}
  (which is taken from~\cite{GGN2015})
  might be confusing,
  since many authors write \qquote{\pre{additive} category}
  to mean a category enriched in abelian groups.
\end{Rem}

\begin{Lem}\cite[Definition~C.1.5.1]{Lurie2018}
  \label{lem:additive-is-in-htpycat}
  Let $\exwaddcat$ be an \inftycat{} with finite products and coproducts.
  Then $\exwaddcat$ is \pre{additive}/additive
  if and only if its homotopy category $\htpycat{\exwaddcat}$
  is \pre{additive}/additive.
\end{Lem}

\begin{Prf}
  The three maps
  defining the \pre{additivity}/additivity of $\exwaddcat$
  in \autoref{def:additive-cat}
  are sent by the functor $\exwaddcat\to\htpycat{\exwaddcat}$
  to the corresponding three maps
  defining the \pre{additivity}/additivity of $\htpycat{\exwaddcat}$.
  Since a map is an equivalence in $\exwaddcat$
  if and only if it is an equivalence
  (\ie isomorphism) in $\htpycat{\exwaddcat}$,
  the result follows.
\end{Prf}

\begin{Lem}
  \label{lem:upper-triangles-are-invertible}
  Let $\exwaddcat$ be an additive \inftycat{}.
  Consider two $n$-tuples
  $(\localx_i)_{i=1}^n$, $(\localy_i)_{i=1}^n$
  of objects of $\exwaddcat$
  and a matrix $F=(f_{i,j}\colon \localx_i\to \localy_j)_{i,j=1}^n$
  of maps between them.
  Assume that
  \begin{itemize}
  \item
    all diagonal entries $f_{i,i}\colon \localx_i\to \localy_i$ are equivalences
  \item
    all entries below the diagonal (\ie $f_{i,j}$ with $i>j$)
    factor through a zero object $\zerobj\in\exwaddcat$.
  \end{itemize}
  Then $F$, viewed as a map
  \begin{equation}
    F\colon
    \coprod_{i=1}^n \localx_i
    \xra{\simeq}
    \prod_{i=1}^n \localy_i,
  \end{equation}
  is an equivalence.
\end{Lem}

\begin{Prf}
  By passing to the homotopy category,
  we may reduce to the case of ordinary additive categories;
  \autoref{lem:upper-triangles-are-invertible} is standard in this case.
\end{Prf}

\subsection{Weakly idempotent complete \inftycats{}}

\newcommand{\localp}{p}
\newcommand{\localx}[1][]{\manyprime[#1]{X}}

Recall that an additive $1$-category $\exabcat$
is called idempotent complete (or Karoubian)
if every idempotent endomorphisms
$\localp\colon \localx\to\localx$
induces a direct sum decomposition
$\localx\cong \Im{\localp}\oplus\Ker{\localp}$.
If $\exabcat$ is embedded as a full additive subcategory
of some abelian category,
this amounts to saying that $\exabcat$ is closed under summands,
in particular, every abelian category is idempotent complete.

In the \inftycategorical{} world,
the situation is a bit less favorable;
for instance, even stable \inftycats{}
are not idempotent complete in general\footnote{
  Splitting a $1$-categorical idempotent $\localp$
  amounts to computing the kernels of $\localp$ and of $\Id-\localp$
  which exist in any abelian category.
  In contrast,
  the splitting a (coherent) \inftycategorical{} idempotent
  must be computed as an \emph{infinite} limit
  which is not always possible.
}.
Fortunately for the purposes of this paper,
the weaker condition of \buzzword{weak} idempotent completeness will suffice.
While idempotent completeness is a way to say
that the category is \qquote{closed under summands},
weakly idempotent completeness should be read as
\qquote{closed under direct complements};
in other words $\exabcat$ is weakly idempotent complete additive
if for each $\localx\in\exabcat$
and each direct sum decomposition $\localx\cong\localx[1]\oplus\localx[2]$
(in some ambient abelian category)
we have
$\localx[1]\in\exabcat$
if and only if
$\localx[2]\in\exabcat$.
One way to intrinsically make this definition without reference
to any ambient category is to say that an additive category $\exabcat$
is weakly idempotent complete
if each retraction (\aka split epimorphism) has a kernel
and each section (\aka split monomorphism) has a cokernel
(see for instance \cite[A.5.1]{ThomasonTrobaugh1990}
and \cite[Definition~7.2]{Buehler2010}).

Next, we define weak idempotent completeness in the \inftycategorical{} setting.
Let $\exptdcat$ be a pointed \inftycat{}.
A \introduce{section-retraction pair} in $\exptdcat$,
is a composable pair $(r,s)$
of maps in $\exptdcat$
whose composite $r\circ s$
is an equivalence.
We say that two section-retraction pairs
$(r,s)$ and $(r',s')$ are \introduce{complementary},
if they fit in a commutative diagram
\begin{equation}
    \label{eq:take-kernel-cokernel}
    \cdcomplsecrec
    \bullet
    \bullet
    \bullet
    \bullet
    \bullet
    {s}{r}{s'}{r'}
\end{equation}
where all squares are biCartesian
(\ie both a pushout and a pullback);
more precisely, we say that the diagram \eqref{eq:take-kernel-cokernel}
exhibits $(r',s')$ as the complement of $(r,s)$, and vice versa.
\begin{Rem}
  It is not hard to show using Kan extensions
  that the evident forgetful functor
  \begin{equation}
    \displayset{diagrams \eqref{eq:take-kernel-cokernel} in $\exptdcat$}
    \xra{\simeq}
    \displayset{section-retraction pairs in $\exptdcat$
      which admit a complement}
  \end{equation}
  is an equivalence \inftycats{}.
  This is the sense in which the complement of a section-retraction pair
  (together with the data exhibiting it as complementary)
  is essentially unique (if it exists).
\end{Rem}

\begin{Def}
  \label{def:weakly-idempotent-complete}
  A pointed \inftycat{} $\exptdcat$
  is called \introduce{weakly idempotent complete}
  if every section-retraction pair
  has a complement.
\end{Def}

\begin{Rem}
  When $\exptdcat=\exabcat$ is an additive $1$-category,
  specifying a diagram \eqref{eq:take-kernel-cokernel}
  amounts to exhibiting $s'$ as the kernel of $r$
  and $r'$ as the cokernel of $s$.
  Hence in this case \autoref{def:weakly-idempotent-complete}
  agrees with the classical notion of weak idempotent completeness.
\end{Rem}

\begin{Expl}
  Every stable \inftycat{} is weakly idempotent complete.
  More generally, each stable \inftycat{}
  gives rise to many examples
  by passing to \sub{categories}
  which are closed under direct complements.
\end{Expl}

\section{The main theorem}

\subsection{\cooltechs{}}

\label{sec:setup}
\newcommand{\localf}{f}

In this section we describe the axiomatic framework
of \cooltechs{}
which encompasses---%
and is essentially equivalent---%
to the setting of \nameLS\cite{LackStreet2015};
see \autoref{sec:comparison-Lack-Street} for a detailed comparison.
Similar ideas were already present in prior work of
\Slominska{}\cite{Slominska2004,Slominska2011}
and of Helmstutler\cite{Helmstutler2014}
(\cf \autoref{rem:comparison-Helmstutler}).

Let $\localB$ be a category
equipped with two \sub{categories} $\epis,\dualepis\subset\localB$,
each of which contains all isomorphism (in particular all objects).
Arrows in $\epis$ and $\dualepis$
are called \introduce{\wepis{}}
and \introduce{\wdualepis{}},
respectively;
we depict them with the symbols $\arepi$
(a two-headed arrow)
and $\ardualepi$ (a tailed arrow),
respectively.
For each $\localb\in\localB$ we denote by
$\episof{\localb}$ the category of \wepis{} under $\localb$\footnote{
  The category $\episof{\localb}$ is nothing but the undercategory
  $\undercat{\localb}{\epis}$
  (where $\localb$ is viewed as an object of $\epis$).
  We do not use the latter notation because it can unfortunately
  be confused with the undercategory
  $
  \undercat{\localb}{\epis}
  =
  \undercat{\localb}{\localB}
  \fiberproduct{\localB}
  \epis
  $
  (where $\localb$ is viewed as an object of $\localB$).
};
similarly, we denote by $\dualepisof{\localb}$ the category
of \wdualepis{} over $\localb$.

We make the following auxiliary definitions:
\begin{itemize}
\item
  We call an arrow in $\localB$ \introduce{singular}
  if it lies in the right ideal
  $\singulars\coloneqq\dualepisneq\circ\localB$ generated
  by the \non{invertible} \wdualepis{}.
\item
  An arrow which is not singular is called \introduce{regular};
  we denote by $\regulars\coloneqq \localB\setminus \singulars$
  the set of regular arrows.
\item
  We call an arrow a \introduce{\wmono{}}
  if it does not lie in
  the left ideal generated by the \non{invertible} \wepis{}.
  We denote by $\monos\coloneqq \localB\setminus (\localB\circ\episneq)$
  the set of \wmonos{}.
\item
  For each $\localb\in\localB$ we have a pairing
  $\blank\circ\blank\colon\episof{\localb}\times\dualepisof{\localb}
  \to \Ar{\localB}$
  given by composition
  (where $\Ar\localB$ denotes the category of arrows in $\localB$).
  We denote by
  \begin{equation}
    \matrixcomponents[\localb]{}{}\colon
    \pi_0\episof{\localb}\times\pi_0\dualepisof{\localb}
    \lra \pi_0\Ar{\localB}
  \end{equation}
  the induced pairing on isomorphism classes.
\end{itemize}

\begin{Def}
  \label{def:cooltech}
  \newcommand\localn{n}
  The datum
  $\excooltech\coloneqq(\localB,\epis,\dualepis)$
  is called
  \begin{itemize}
  \item
    A \introduce{\cooltech{}}\footnote{
    Unsurprisingly, DK stands for Dold--Kan.}
    if it satisfies the following properties
    (using the auxiliary notation introduced above):
    \begin{enumerate}[label=\axiomlabelstyle{T}, ref=\axiomrefstyle{T}]
    \item
      \label{item:cooltech-fact-property}
      Every arrow $\localf$ of $\localB$ can be written uniquely
      (up to unique isomorphism)
      as a composition of the form
      \begin{equation}
        \label{eq:three-fold-factorization}
        \begin{tikzcd}[column sep = large]
          \bullet\ar[r,"{\exepi[1]\in\epis}",twoheadrightarrow]
          &\bullet\ar[rr,"{\ol{\localf}\in(\mcapr)}"]
          &
          &\bullet\ar[r,"{\exdualepi\in\dualepis}",rightarrowtail]
          &\bullet
        \end{tikzcd}
      \end{equation}
    \item
      \label{item:cooltech-UT-matrix}
      For each $\localb\in\localB$,
      the pairing $\matrixcomponents[\localb]{}{}$
      can be described by a finite square matrix which is
      \roughly{unipotent upper triangular modulo \non{isomorphisms}},\\
      \ie there is a number $\localn\geq 1$ and bijections
      $\pi_0\episof{\localb}
      \cong
      \set{1,\dots,\localn}
      \cong
      \pi_0\dualepisof{\localb}$,
      such that the pairing
      $\matrixcomponents[\localb]{}{}$
      induces an $\localn\times \localn$-matrix
      \newcommand\localph{?}
      \begin{equation}
        \left(
          \begin{matrix}
            \simeq & \localph& \cdots&\localph\\
            \not\simeq& \ddots & \ddots &\vdots\\
             \vdots&\ddots&\ddots & \localph\\
            \not\simeq&\cdots &\not\simeq&\simeq \\
          \end{matrix}
        \right)
      \end{equation}
      with values in $\pi_0\Ar{\localB}$ which has
      invertible arrows on the diagonal and
      \non{invertible} arrows below the diagonal
      (there is no condition on the arrows above the diagonal).
    \item
      \label{item:cooltech-epis-dual-epis-cat}
      The set
      $
      \dualepis\circ\epis
      $
      is closed under composition.
    \item
      \label{item:cooltech-monos-cat}
      The composition of two regular \wmonos{}
      is a (not necessarily regular) \wmono{},\\
      \ie $(\mcapr)\circ(\mcapr)\subset\monos$
    \item
      \label{item:cooltech-sing-ideal}
      The singular arrows form a left module over $\monos$,
      \ie we have
      $\monos\circ\singulars \subseteq \singulars$.
    \end{enumerate}
  \item
    a \introduce{\reallycooltech{}}
    if it if satisfies all axioms
    \ref{item:cooltech-fact-property}--\ref{item:cooltech-sing-ideal}
    above and the matrix in
    \ref{item:cooltech-UT-matrix}
    can even be made \emph{diagonal} modulo \non{isomorphisms}.
  \item
    \introduce{reduced}
    if $\localB=\dualepis\circ\epis$.
    \qedhere
  \end{itemize}
\end{Def}

The following observations follow immediately from \autoref{def:cooltech}.

\begin{Lem}
  \label{lem:cooltech-basics}
  Let $\excooltech=(\localB,\epis,\dualepis)$ be a \cooltech{}.
  \begin{enumerate}
  \item
    \label{item:dualepi-epi-section-retraction}
    For each $\localb$ there is a unique bijection
    $(\blank)^\dual\colon\pi_0\episof{\localb}\llra\pi_0\dualepisof{\localb}$
    such that for each $\exepi\in \episof{\localb}$
    the composition $\exepi\circ \exdualepi$
    is an isomorphism in $\localB$.
  \item
    Every \wepi{} is a split epimorphism
    and every \wdualepi{} is a split monomorphism in $\localB$.
  \item
    For each $\localb\in\localB$,
    the categories $\episof{\localb}$ and $\dualepisof{\localb}$
    are both (equivalent to) posets.
  \item
    Both $\regulars$ and $\monos$ contain all isomorphisms of $\localB$.
  \item
    An arrow $\localB$
    decomposed as in \eqref{eq:three-fold-factorization}
    is regular if and only if
    the component $\exdualepi\in\dualepis$ is invertible
    and it is a \wmono{} if and only if
    the component $\exepi[1]\in\epis$ is invertible.
  \item
    We have $\monos=(\mcapr)\circ\epis$
    and $\regulars=\dualepis\circ(\mcapr)$.
  \item
    The datum
    $\ol{\excooltech}\coloneqq(\ol{\localB},\epis,\dualepis)$
    is again a \cooltech{} which is automatically reduced.
  \item
    If $\excooltech$ is reduced then we have
    $\monos=\dualepis$
    and
    $\regulars=\epis$
    and
    $\mcapr=\grpdcore{\localB}$.
  \item
    If $\excooltech$ is reduced
    then the dual datum
    $\excooltech^\op\coloneqq(\localB^\op,(\dualepis)^\op,\epis^\op)$
    is again a (reduced) \cooltech{}.
    \qedhere
  \end{enumerate}
\end{Lem}

\begin{Prf}
  Straightforward and left to the reader.
\end{Prf}

Each \cooltech{}
$\excooltech=(\localB,\epis,\dualepis)$
induces a canonical partial order
$\leq$ on the set $\pi_0\localB$
by declaring $\localb[1]\leq\localb$
if there exists
a \wdualepi{} $\localb[1]\ardualepi\localb$
or equivalently (by \ref{item:dualepi-epi-section-retraction})
an \wepi{} $\localb\arepi\localb[1]$.
To see that $\leq$ is antisymmetric
(\ie
 $\localb\leq\localb[1]\leq\localb$
implies $\localb\cong\localb[1]$)
choose an \wepi{}
$\exepi\colon\localb[1]\arepi\localb$
and an \wepi{}
$\localb\arepi\localb[1]$:
the induced maps
$\blank\circ\exepi\colon\pi_0\episof{\localb}\hra\pi_0\episof{\localb[1]}$
and
$\pi_0\episof{\localb[1]}\hra\pi_0\episof{\localb}$
are injective
because \wepis{} are (split) epimorphisms.
Since the sets
$\pi_0\episof{\localb}$ and $\pi_0\episof{\localb[1]}$
are finite by \ref{item:cooltech-UT-matrix},
this implies that $\exepi\circ\blank$ is a bijection;
hence $\exepi$ is a split monomorphism
because $\Id_{\localb[1]}$ lies in the image of $\exepi\circ\blank$;
hence $\exepi$ is an isomorphism.

For each $\localb\in\localB$ the set
$\setP{\localb[1]\in\pi_0\localB}{\localb[1]\leq\localb}$
of predecessors of $\localb$ is finite by
\ref{item:cooltech-UT-matrix},
hence the poset $(\pi_0\localB,\leq)$ is suited for inductive arguments.

\subsection{Key constructions}

\label{sec:key-constructions}
\begin{Cstr}
  \label{cstr:quotient-A}
  Assume that $\excooltech$ is a \cooltech{}.
  We define a pointed category
  $\localAp=\localApof{\excooltech}$ as the quotient
  \begin{equation}
    \label{eq:def-A}
    \localAp\coloneqq\quotcat{\monos}{\monos\cap\singulars}
  \end{equation}
  of $\monos$ by the two-sided ideal $\monos\cap\singulars$.
  Explicitly:
  \begin{itemize}
  \item
    The pointed category $\localAp$
    has a zero object $\zerobj$
    and for each object $\localb\in\localB$ an object $\localbA\in\localA$
  \item
    For every pair of objects $\localbA[1],\localbA\in\localA$,
    we have the hom-set
    \begin{equation}
      \localAp(\localbA[1],\localbA)
      \coloneqq
      \frac{\monos(\localb[1],\localb)}{(\monos\cap\singulars)}
      =
      (\mcapr)(\localb[1],\localb)
      \disjunion
      \set{\localbA[1]\to\zerobj\to\localbA}.
    \end{equation}
  \item
    Composition in $\localAp$ is induced by composition in $\localB$;
    it is \well{defined} because of
    \ref{item:cooltech-monos-cat}
    and
    \ref{item:cooltech-sing-ideal}.
  \end{itemize}
  For convenience we write $\localA$
  for the full \sub{category} of $\localAp$
  spanned by all objects except the zero object $\zerobj$.
\end{Cstr}

\begin{Rem}
  \label{rem:cstr-when-Ap-is-free}
  A particularly simple case of \autoref{cstr:quotient-A}
  occurs when the set $\mcapr$ of regular \wmonos{}
  is closed under composition.
  In this case,
  $\mcapr$ is a \sub{category} of $\localB$
  and the quotient
  $\localAp\coloneqq
  \quotcat{\monos}{\monos\cap\singulars}
  \cong
  \quotcat{\mcapr}{\emptyset}
  \cong
  \quotcat{\localA}{\zeroideal}
  $
  is simply the free pointed category
  on the category $\mcapr$.
\end{Rem}

\begin{Not}
  To minimize the potential confusion, we adopt the following conventions:
    Objects in $\localA$ are denoted by
    $\locala,\locala[1],\locala[2]$.
    Objects in $\localB$ are denoted by
    $\localb,\localb[1],\localb[2]$.
    Given an object $\locala\in\localA$,
    we denote by $\localaB$ the corresponding object in $\localB$.
    \qedhere
\end{Not}

We now come to the key construction of this article.

\begin{Cstr}
  \label{cstr:wrapcat}
  Let $\excooltech$ be a \cooltech{}.
  We define the pointed category
  $\wrapcat=\wrapcatof{\excooltech}$
  as the \roughly{upper triangular} category
  \begin{equation}
    \wrapcat\coloneqq
    \left(
      \begin{matrix}
        \localAp& \regularsP\\
        \zerobj & \localBp
      \end{matrix}
    \right)
    \coloneqq
    \left(
      \begin{matrix}
        \frac{\monos}{\singulars}& \lquot{\localB}{\singulars}\\
        \zerobj & \localBp
      \end{matrix}
    \right)
  \end{equation}
  associated to the
  $\localAp$-$\localBp$-bimodule
  $\regularsP\coloneqq\lquot{\localB}{\singulars}$.
  More precisely,
  the category $\wrapcat$ is given explicitly as follows:
  \begin{itemize}
  \item
    The objects of $\wrapcat$ are given by
    the objects $\locala\in\localA$,
    the objects $\localb\in\localB$
    and a zero object $\zerobj$;
    in other words we have
    $
    \Ob{\wrapcat}
    \coloneqq
    \Ob{\localAp}\pushout{\set{\zerobj}}\Ob{\localBp}
    $.
  \item
    The hom-sets in $\wrapcat$ between two objects of $\localAp$
    or between two objects of $\localBp$
    are inherited from $\localAp$
    or from $\localBp$,
    respectively.
  \item
    The only arrow in $\wrapcat$ from an object
    $\locala\in\localAp$
    to an objects
    $\localb\in\localBp$
    is the zero arrow
    $\locala\to\zerobj\to\localb$
  \item
    The set of arrows in $\wrapcat$ from
    $\localb\in\localB$ to $\locala\in\localA$
    is defined to be
    \begin{equation}
      \wrapcat(\localb,\locala)
      \coloneqq
      \regularsP(\localb,\locala)
      \coloneqq
      \lquot{\localB(\localb,\localaB)}{\singulars}
      =
      \regulars(\localb,\localaB)
      \disjunion
      \set{\localb\to\zerobj\to\locala}
    \end{equation}
  \item
    Composition in $\wrapcat$ is induced by the composition
    in $\localAp$ and in $\localBp$;
    the composition
    \begin{equation}
      \localAp(\locala,\locala[1])
      \times
      \regularsP(\localb,\locala)
      \times
      \localBp(\localb[1],\localb)
      \lra
      \regularsP(\localb[1],\locala[1])
    \end{equation}
    is \well{defined} because
    $\monos\circ\singulars\circ\localB\subseteq\singulars$.
  \end{itemize}
  The pointed category $\wrapcat$ comes equipped
  with the two fully faithful embeddings
  \begin{equation}
    \localBp\lhra\wrapcat\lhla\localAp;
  \end{equation}
  for convenience we identify $\localBp$ and $\localAp$
  with their images in $\wrapcat$.
\end{Cstr}

\begin{Not}
  We denote by
  $\canmapBA\colon\localaB\to\locala$
  the arrow corresponding to the identity
  $\Id_{\localaB}\in\regulars(\localaB,\localaB)$.
  For every \non{zero} arrow
  $\localu\colon\localb\to\locala$ in $\wrapcat$
  we denote by $\localuB\in\regulars(\localb,\localaB)$
  the corresponding regular arrow in $\localb$;
  in other words,
  $\localuB\colon\localb\to\localaB$ is the unique arrow
  satisfying $\canmapBA{\localuB}=\localu$.
\end{Not}

\begin{Rem}
  Assumptions
  \ref{item:cooltech-monos-cat} and \ref{item:cooltech-sing-ideal}
  are needed to guarantee that
  \autoref{cstr:quotient-A} and \autoref{cstr:wrapcat}
  are \well{defined}.
  In many examples
  $\monos$ is actually a \sub{category} of $\localB$;
  in this case $\monos\cap\singulars$
  is a two-sided ideal in $\monos$ in the usual sense
  and \autoref{cstr:quotient-A}
  becomes an instance of
  \autoref{cstr:quotient-pointed-cat}.
  The notation in
  \autoref{cstr:quotient-A} and \autoref{cstr:wrapcat}
  should be understood with this
  more special (but still very general) case in mind.
\end{Rem}

\subsection{Statement}

We now state the main theorem of this article.

\begin{Thm}[Homotopy coherent correspondences of Dold--Kan type]
  \label{thm:main-thm}
  Let $\excooltech=(\localB,\epis,\dualepis)$ be a \cooltech{}
  with associated pointed category
  $\localAp=\localAp(\excooltech)$.
  \begin{enumerate}[label=(\alph*)]
  \item
    \label{it:main-Kan-extensions-exist}
    For any \waddjectives{} \inftycat{}
    $\exwaddcat$,
    the restriction functors
    \begin{equation}
      \Res\colon\waddiagramsPof{\wrapcat}
      \lra\waddiagramsPof{\localBp}
      \intxt{and}
      \Res\colon\waddiagramsPof{\wrapcat}
      \lra\waddiagramsPof{\localAp}
    \end{equation}
    from \autoref{cstr:wrapcat}
    admit a left adjoint $\LKE$ (left Kan extension)
    and a right adjoint $\RKE$ (right Kan extension), respectively.
  \item
    \label{it:main-master-equivalence}
    The composite adjunction
    \begin{equation}
      \label{eq:master-span-LRKE}
      \begin{tikzcd}
        \waddiagramsPof{\localBp}
        \tikzcdadjunction{\LKE}{\Res}
        &
        \waddiagramsPof{\wrapcat}
        \tikzcdadjunction{\Res}{\RKE}
        &
        \waddiagramsPof{\localAp}
        \\
      \end{tikzcd}
    \end{equation}
    is an adjoint equivalence of \inftycats{}.
  \item
    \label{it:main-naturality}
    The adjoint equivalence
    \eqref{eq:master-span-LRKE}
    is natural in $\exwaddcat$ with respect to additive functors.
  \item
    \label{it:main-colim-lim-formula}
    Consider a pointed functor
    $\exptdiagramB\colon\localBp\to\exwaddcat$
    and denote by
    $\exptdiagramA\colon\localAp\to\exwaddcat$
    the pointed functor corresponding to $\exptdiagramB$
    under the equivalence \eqref{eq:master-span-LRKE}.
    Then for each $\locala\in\localA$ the canonical maps
    \begin{equation}
      \label{eq:main-colim-lim-sec-rec}
    {
      \colim\limits_{\localb\in\dualepisneqof{\localaB}}
      \localevalB{\exdiagramB}{\localb}
    }
    \lra
    \localevalB{\exdiagramB}{\localaB}
    \lra
    {
      \lim\limits_{\localb\in\episneqof{\localaB}}
      \localevalB{\exdiagramB}{\localb}
    }
    \end{equation}
    form a section-retraction pair
    with complement equivalent to
    $\localevalA{\exptdiagramA}{\locala}$.
    \qedhere
  \end{enumerate}
\end{Thm}

\begin{Rem}
  \label{rem:main-thm-dualizes}
  The notions of \ppre{additivity} and weak idempotent completeness
  are manifestly \self{dual}.
  Therefore in \autoref{thm:main-thm} (and all of the results below)
  we can replace the target \inftycat{} by its opposite,
  or, equivalently,
  $\localBp$ by $(\localBp)^\op$
  and
  $\localAp$ by $(\localAp)^\op$.
\end{Rem}

\begin{Rem}
  Since we are not assuming that our target category $\exwaddcat$
  has finite limits or colimits,
  it is not clear \latin{a priori}
  that the limits/colimits indicated in
  \eqref{eq:main-colim-lim-sec-rec}
  even exist;
  part of the statement of
  \autoref{thm:main-thm}~\ref{it:main-colim-lim-formula}
  is that they do.
  Similarly,
  \ref{it:main-Kan-extensions-exist}
  is not automatic;
  in fact, the heart of the proof of \autoref{thm:main-thm}
  is an explicit inductive pointwise construction
  of the Kan extensions \eqref{eq:master-span-LRKE}
  in the case where $\excooltech$ is reduced
  (see \autoref{prop:two-sided-Kan-extension}).
\end{Rem}

\begin{Cor}
  \label{cor:main-thm-without-pt}
  In the situation of \autoref{thm:main-thm},
  the span
  $\localB\subset\localBp\hra\wrapcat\hla\localAp$
  induces a natural equivalence
  \begin{equation}
    \label{eq:equiv-DB-BpAp}
    \waddiagramsof{\localB}\xlra{\simeq}\waddiagramsPof{\localAp}
  \end{equation}
  of \inftycats{}
  for each \waddjectives{} \inftycat{} $\exwaddcat$.
\end{Cor}

\begin{Prf}
  Compose the equivalence of \autoref{thm:main-thm}
  with the natural equivalence
  \begin{equation}
    \waddiagramsof{\localB}\xla{\simeq}\waddiagramsPof{\localBp}
  \end{equation}
  produced by the universal property
  of the free pointed category $\localB\to\localBp$.
\end{Prf}

\begin{Rem}
  \label{rem:thm-when-Ap-is-free}
  In the situation of \autoref{rem:cstr-when-Ap-is-free},
  where $\localAp=\addbasept{(\mcapr)}$
  is a \emph{free} pointed category,
  we can simplify the statement of \autoref{cor:main-thm-without-pt}
  even more and obtain a natural equivalence
  \begin{equation}
    \waddiagramsof{\localB}
    \xlra{\simeq}
    \waddiagramsof{\mcapr}
  \end{equation}
  between ordinary (\ie \non{pointed}) \inftycats{} of diagrams.
  All equivalences discussed in \autoref{sec:cooltechs-from-factorizations}
  are of this form.
\end{Rem}

Specializing \autoref{cor:main-thm-without-pt}
to the $1$-categorical case,
we recover the main theorem of \nameLS{}.

\begin{dCor}\cite[Theorem~6.8]{LackStreet2015}
  \label{cor:main-thm-1-categorical}
  Each \cooltech{} $\excooltech=(\localB,\epis,\dualepis)$
  induces a natural equivalence
  \begin{equation}
    \label{eq:equiv-DB-BpAp}
    \Fun(\localB,\exabcat)
    \xlra{\simeq}
    \ptdFun{\localApof{\excooltech}}{\exabcat}
  \end{equation}
  of categories
  for each weakly\footnote{
    To be precise, \nameLS{} assume
    $\exabcat$ to be idempotent complete.
  } idempotent complete additive category $\exabcat$.
\end{dCor}

\begin{Rem}
  \label{rem:naturality-to-htpy-cat}
  Since the functor
  $\exwaddcat\to\htpycat\exwaddcat$
  to the homotopy category is additive,
  the naturality of equivalence
  \eqref{eq:equiv-DB-BpAp}
  implies the existence of a commutative square
  \begin{equation}
    \cdsquareOpt
    {\waddiagramsof{\localB}}
    {\waddiagramsPof{\localAp}}
    {\Fun(\localB,\htpycat\exwaddcat)}
    {\ptdFun{\localAp}{\htpycat\exwaddcat}}
    {"\simeq"}
    {}{}
    {"\simeq"}
  \end{equation}
  where the lower equivalence
  is an instance of \autoref{cor:main-thm-1-categorical}.
\end{Rem}

\begin{Rem}
  \label{rem:thm-when-cooltech-is-really}
  If the \cooltech{} $\excooltech$ is \reallycool{},
  then in all of the results above
  one can weaken the assumption on $\exwaddcat$
  and only require it to be weakly idempotent complete and \emph{pre}additive.
  Indeed, the additivity of $\exwaddcat$
  is only used once (in the proof of \autoref{prop:two-sided-Kan-extension})
  to invert certain upper triangular matrices in $\exwaddcat$
  obtained from the matrices
  $\matrixcomponents[\localb]{}{}$
  defined in \ref{item:cooltech-UT-matrix};
  if $\excooltech$ is \reallycool{}
  then these matrices in $\exwaddcat$ are diagonal,
  hence inverting them only requires \pre{additivity}.
  See also \autoref{rem:when-matrix-in-prop-is-diagonal}.
\end{Rem}

\begin{Rem}
  \label{rem:tot-fib-cofib-formulas}
  \autoref{thm:main-thm}~\ref{it:main-colim-lim-formula}
  implies that one can compute the value
  $\localevalA\exptdiagramA{\locala}$
  of the diagram
  $\exptdiagramA\colon\localAp\to\exwaddcat$
  at an object $\locala\in\localA$
  in two seemingly unrelated ways:
  as a total fiber
  \begin{equation}
    \label{eq:tot-fib-formula}
    \localevalA\exptdiagramA{\locala}
    \simeq
    \totfib
    \left(
      \episof{\localaB}\to\localB\xra{\exdiagramB}\exwaddcat
    \right)
  \end{equation}
  along the \wepis{},
  or as a total cofiber
  \begin{equation}
    \label{eq:tot-cofib-formula}
    \localevalA\exptdiagramA{\locala}
    \simeq
    \totcof
    \left(
      \dualepisof{\localaB}\to\localB\xra{\exdiagramB}\exwaddcat
    \right).
  \end{equation}
  along the \wdualepis{}.
\end{Rem}

\section{Examples}

\subsection{The \inftycategorical{} \DKcorr{}}

\label{sec:example-Delta}
\newcommand{\numnnCh}[1]{\ul{#1}}
\newcommand{\localf}{f}
\newcommand{\localg}{g}
\newcommand{\locali}{i}
\newcommand{\Dcoolmin}{
  \renewcommand{\dualepis}{\epis^\dual_{\min}}
  \renewcommand{\dualepisneq}{\Error}
  \renewcommand{\localB}{\Delta}
  \renewcommand{\excooltech}{\BB^\Delta_{\min}}
}

\newcommand{\Dcoolmax}{
  \renewcommand{\dualepis}{\epis^\dual_{\max}}
  \renewcommand{\dualepisneq}{\Error}
  \renewcommand{\localB}{\Delta}
  \renewcommand{\excooltech}{\BB^\Delta_{\max}}
}

\newcommand{\DKmin}{\mathrm{DK}^{\min}}
\newcommand{\DKmax}{\mathrm{DK}^{\max}}
\newcommand{\caninvolution}{\updownarrow}

We explain how to equip the simplex category
$\localB=\Delta$ with the structure of a \cooltech{};
\cf \cite[Example 3.2]{LackStreet2015}.
Recall that $\Delta$ is the category of
finite \non{empty} linearly ordered sets and weakly monotone maps between them.
We denote by $\numD{\nn}$ the standard ordinal $\set{0<1<\dots<\nn}$;
every object of $\Delta$ is of this form up to unique isomorphism.
Let $\epis\subset\Delta$
be the wide subcategory of surjective maps
and let $\Dcoolmin\dualepis\subset\Delta$
be the wide subcategory of those injectives maps that preserve minimal elements.
The following observations are straightforward to verify
and imply that $\Dcoolmin\excooltechdef$
is a \cooltech{}:
\begin{itemize}
  \Dcoolmin
\item
  A map $\localf\colon\numD{\nn}\to\numD{\mm}$ is singular
  if and only if there is a \non{minimal} element
  of $\numD{\mm}$ which is not in the image of $\localf$.
\item
  The set $\monos$ of \wmonos{} consists precisely
  of the injective maps in $\Delta$.
  Since $\monos$ is closed under composition,
  \ref{item:cooltech-monos-cat} is satisfied.
\item
  The set $\mcapr$ of regular monos consists of the identities
  and the $0$-th coface maps
  \begin{equation}
    \cofacemap{0}
    \colon
    \numD{\nn-1}\cong\set{1<\dots<\nn}
    \lhra\numD{\nn}.
  \end{equation}
  Note that $\mcapr$ is not closed under composition.
\item[\ref{item:cooltech-fact-property}]
  Each map $\numD{\nn}\to\numD{\mm}$ in $\Delta$
  admits a unique (up to unique isomorphism) factorization of type
  $\dualepis\circ(\mcapr)\circ\epis$,
  namely
  \begin{equation}
    \numD{\nn}
    \lthra{}
    (\Im{\localf})
    \lhra
    \left(\set{0}\cup\Im{\localf}\right)
    \rightarrowtail
    \numD{\nn}.
  \end{equation}
\item[\ref{item:cooltech-epis-dual-epis-cat}]
  The set $\dualepis\circ\epis$
  consists of the minimum-perserving arrows in $\Delta$,
  hence $\dualepis\circ\epis$ is closed under composition.
\item[\ref{item:cooltech-sing-ideal}]
  If $0\neq\locali\in\numD{\mm}$ is a \non{minimal} element
  which is not in the image of
  $\localf\colon\numD{\nn}\to\numD{\mm}$
  then, for each injective map
  $\localg\colon\numD{\mm}\to\numD{\mm'}$,
  the element $0\neq\localg(\locali)\in\numD{\mm'}$
  is not minimal and not contained in the image of $\localg\circ\localf$.
\item[\ref{item:cooltech-UT-matrix}]
  For each $\numD{\nn},\numD{\mm}\in\Delta$
  we have a bijection
  \begin{equation}
    (\blank)^\dual
    \colon
    \epis(\numD\nn,\numD\mm)
    \xra{\cong}
    \dualepis(\numD\mm,\numD\nn)
  \end{equation}
  which sends a surjection $\exepi\colon\numD\nn\to\numD\mm$
  to its minimal section
  $
  \exdualepi\colon\numD{\nn}\to\numD\nn
  $
  given by the formula
  $
  \locali\mapsto\min\preimage{\exepi}{\locali}.
  $
  A composition
  $\exepi[1]\circ\exdualepi$
  of an \wepi{}
  $\exepi[1]\colon\numD\nn\arepi\numD\mm$
  with a \wdualepi{}
  $\exdualepi\colon\numD\nn\ardualepi\numD\mm$
  is
  \begin{itemize}
  \item
    an isomorphism if $\exepi[1]=\exepi$
  \item
    not an isomorphism if
    $\exdualepi\not\geq\exdualepi[1]$
    poinwise as maps $\numD\mm\to\numD\nn$
  \end{itemize}
  (note that we make no claim when $\exdualepi[1]<\exdualepi$).
  Hence, for each $\numD{\nn}\in\Delta$,
  the lexicographic ordering on $\pi_0\dualepisof{\numD\nn}$
  makes the matrix
  \begin{equation}
    \pi_0\dualepisof{\numD\nn}
    \times
    \pi_0\dualepisof{\numD\nn}
    \xla[(\blank)^\dual]{\cong}
    \pi_0\episof{\numD\nn}
    \times
    \pi_0\dualepisof{\numD\nn}
    \xra{\matrixcomponents[\numD\nn]{}{}}
    \pi_0\Ar\Delta
  \end{equation}
  into upper triangular shape modulo \non{isomorphisms}.
\end{itemize}

\begin{Expl}
  \label{ex:example-matrix-Delta}
  \newcommand{\grouping}[1]{(#1)}
  The matrix
  $\matrixcomponents[\numD{2}]{}{}$,
  can be depicted as follows
  \begin{equation}
    \newcommand{\localhl}[1]{\textcolor{red}{\mathbf{#1}}}
    \begin{array}[c]{c|cccc}
      &0&01&02&012\\
      \hline
      \grouping{012}
      &\localhl{0}&\grouping{01}&\grouping{12}&\grouping{012}\\
      0\grouping{12}
      &0&\localhl{01}&\localhl{02}&0\grouping{12}\\
      \grouping{01}2
      &0&\grouping{01}&\localhl{02}&\grouping{01}2\\
      012
      &0&01&02&\localhl{012}
    \end{array}
  \end{equation}
  where the rows are labeled by equivalence classes of \wepis{}
  $\numD{2}\arepi\numD{\mm}$
  (written by grouping elements with the same image);
  dually, the rows are labeled by equivalence classes of \wdualepis{}
  $\numD{\mm}\arepi\numD{2}$
  (written by listing the elements in the image).
  The isomorphisms are highlighted in red/bold,
  showing that the matrix is---modulo \non{isomorphisms}---%
  unipotent upper triangular but not diagonal.
  In particular, this example shows that
  the \cooltech{} $\Dcoolmin\excooltech$ is not \really{}.
\end{Expl}

There is an equivalence
$\nnCh\coloneqq
\frac{\BN}{\idealtwodiffop}
\xra{\simeq} \localAp\coloneqq \frac{\monos}{\mcapr}$
of pointed categories
which is given on objects by
$\numnnCh{\nn}\mapsto \numD{\nn}$
and is determined on morphisms by sending the arrow
$\numnnCh{\nn}\to\numnnCh{\nn+1}$
to the $0$-th coface map
$\cofacemap{0}\colon\numD{\nn}\to\numD{\nn+1}$.
Applying \autoref{cor:main-thm-without-pt}
to the \cooltech{} $\Dcoolmin\excooltech$
thus establishes a natural equivalence
\begin{equation}
  \Fun(\Delta,\exwaddcat)
  \xlra{\simeq}
  \ptdFun{\nnCh}{\exwaddcat}
  =
  \nnChof{\exwaddcat}
\end{equation}
for each \waddjectives{} \inftycat{} $\exwaddcat$.
Replacing $\exwaddcat$ by its opposite
(which is again \waddjectives{})
yields the more familiar form of
the following \introduce{\inftycategorical{} \DKcorr{}}:

\begin{dCor}
  \label{cor:infty-categorical-DK}
  The \cooltech{}
  $\Dcoolmin\excooltechdef$
  induces a natural equivalence of \inftycats{}
  \begin{equation}
    \label{eq:infty-categorical-DK}
    \Fun(\Dop,\exwaddcat)
    \xlra{\simeq}
    \connChof{\exwaddcat}
  \end{equation}
  between simplicial objects and connective chain complexes in any
  \waddjectives{} \inftycat{} $\exwaddcat$.
\end{dCor}

{
  \Dcoolmax
  The simplex category $\Delta$ is part of a second \cooltech{}
  $\excooltech$,
  where $\epis$ is again the set of surjections
  and $\dualepis$ is the set of maximum-preserving injections in $\Delta$.
}
The \cooltechs{} $\Dcoolmin\excooltech$ and $\Dcoolmax\excooltech$
have canonically isomorphic quotient categories
$
\localAp({\Dcoolmin\excooltech})
\cong
\nnCh
\cong
\localAp({\Dcoolmax\excooltech})
$
and correspond to each other under the canonical involution
$\caninvolution\colon\Delta\xlra{\cong}\Delta$
which reverses the linear order on each object of $\Delta$.
Hence we have a commutative diagram
\begin{equation}
  \cdtriangleoverOpt
  {\waddiagramsof{\Dop}}
  {\waddiagramsof{\Dop}}
  {\connChof{\exwaddcat}}
  {"\caninvolution"',"\simeq",leftrightarrow}
  {"\DKmax","\simeq"',leftrightarrow}
  {"\DKmin"',"\simeq",leftrightarrow}
\end{equation}
which intertwines the corresponding two versions
of the Dold--Kan functor \eqref{eq:infty-categorical-DK}.
\begin{Lem}
  For every additive \inftycat{} $\exwaddcat$,
  the autoequivalence
  $
  \caninvolution\colon
  \waddiagramsof{\Dop}
  \xlra{\simeq}
  \waddiagramsof{\Dop}
  $
  is equivalent to the identity;
  in other words,
  the two Dold--Kan functors
  $\DKmin$ and $\DKmax$ agree
  (up to equivalence).
\end{Lem}

In the $1$-categorical context,
one can check explicitly that
$\DKmin$ and $\DKmax$ both
agree (up to natural isomorphism) with the normalized chain functor,
hence with each other.
For \inftycats{} we provide the following alternative argument:

\begin{Prf}
  Choose a stable \inftycat{} $\exwaddcat\subseteq\exstablecat$
  into which $\exwaddcat$ is embedded as a full additive subcategory;
  for instance the Yoneda embedding
  into
  the \inftycat{} of additive spectral presheaves
  $\exwaddcat^\op\to\Spectra$
  (for more details, see for instance Appendix~C.1.5 in \cite{Lurie2018}).
  Since the involution
  $\caninvolution\colon\Delta\xlra{\cong}\Delta$
  preserves the filtration
  \begin{equation}
    \Deltaleq{0}\subset\Deltaleq{1}
    \subset\cdots\subset
    \Deltaleq{\nn}
    \subset\cdots\subset\Delta,
  \end{equation}
  the functor
  $
  \Fun(\Dop,\exstablecat)
  \to
  \Fun(\BN,\exstablecat),
  $
  which sends a simplicial objects to its sequence of partial colimits,
  is $\caninvolution$-invariant.
  Since $\exstablecat$ is stable,
  Lurie's stable \DKcorr{}
  (which we also discuss in \autoref{sec:comparison-Lurie})
  states that this functor is an equivalence;
  hence the ($\caninvolution$-invariant) composition
  \begin{equation}
    \Fun(\Dop,\exwaddcat)
    \hra
    \Fun(\BN,\exstablecat),
  \end{equation}
  is fully faithful.
  The result follows.
\end{Prf}

\subsection{Categories of partial maps}

\label{sec:cooltechs-from-factorizations}
\newcommand{\exfactm}[1][]{\manyprime[#1]{m}}
\newcommand{\exfactE}{e}
\newcommand{\exfactM}{m}

An important class of examples of
\reallycooltechs{}
arises by considering partial maps
with respect to certain factorization systems;
we sketch here the corresponding discussion by
\nameLS{}\cite[Example 3.1]{LackStreet2015}.
Let $(\factE,\factM)$ be a factorization system
on a category $\factcat$,
\ie
\begin{itemize}
\item
  $\factE$ and $\factM$ are \sub{categories} of $\factcat$
  containing all isomorphisms and
\item
  arrows in $\factcat$ factor,
  uniquely up to unique isomorphism,
  as compositions of the type $\factM\circ\factE$;
\end{itemize}
assume furthermore that
\begin{itemize}
\item
  the arrows in $\factM$ are monomorphisms,
\item
  each object of $\factcat$ has only finitely many $\factM$-\sub{objects} and
\item
  the pullback of an arrow in $\factM$
  along an arbitrary map in $\factcat$ exists
  and lies again in $\factM$.
\end{itemize}
The category $\partialmapsop{\factcat}=\partialmapsop[\factM]{\factcat}$
of \introduce{co-$\factM$-partial maps} is defined to have
\begin{itemize}
\item
  the same objects as $\factcat$;
\item
  morphism in $\partialmapsop{\factcat}$
  are equivalence classes of spans in $\factcat$ of the type
  $\bullet\xla{\factcat}\bullet\xra{\factM}\bullet$,
  \ie where the second leg is required to lie in $\factM$;
\item
  composition in $\partialmapsop{\factcat}$ is that of spans,
  \ie by pullback.
\end{itemize}
We define two wide \sub{categories}
$\epis,\dualepis\subset\partialmapsop{\factcat}$
to consist of the spans of the type
\begin{equation}
  \bullet\xla{\factM}\bullet\xra{\cong}\bullet
  \intxt{and}
  \bullet\xla{\cong}\bullet\xra{\factM}\bullet,
\end{equation}
respectively.
With the notation of \autoref{sec:setup},
the \wmonos{} are the spans of type
$\bullet\xla{\factE}\bullet\xra{\factM}\bullet$.
The set $\regulars$ of regular morphisms consists
of those morphisms in $\partialmapsop{\factcat}$
which are totally defined,
\ie the spans of type $\bullet\xla{\factcat}\bullet\xra{\cong}\bullet$;
hence we have $\regulars\cong\factcat^\op$.

\begin{Lem}
  \label{lem:facotorizations-give-cooltechs}
  The datum
  $\excooltech\coloneqq(\partialmapsop{\factcat},\epis,\dualepis)$
  is a \reallycool{} \cooltech{}.
\end{Lem}

\begin{Prf}
  The proof is straightforward and left to the reader.
\end{Prf}
The regular \wmonos{} in
$\partialmapsop{\factcat}$
are the spans of the type
$\bullet\xla{\factE}\bullet\xra{\cong}\bullet$;
they form a \sub{category} equivalent to $\factE^\op$.
Hence \autoref{rem:cstr-when-Ap-is-free} says
that the pointed category
$\localAp(\excooltech)$
constructed in \autoref{cstr:quotient-A}
is just the free pointed category $\addbasept{\factE^\op}$
on $\mcapr\cong\factE^\op$.

The upshot of this discussion is the following corollary of
\autoref{cor:main-thm-without-pt},
taking into account that
$\localAp=\addbasept{\factE^\op}$ is a free pointed category
(see \autoref{rem:thm-when-Ap-is-free})
and that
the \cooltech{} $\excooltech$ is \really{}
(see \autoref{rem:thm-when-cooltech-is-really}).

\begin{Cor}
  \label{cor:main-thm-for-partial-maps}
  Let $\factcat$ and $(\factE,\factM)$ be as above.
  The \cooltech{}
  $(\partialmapsop[\factM]{\factcat},\epis,\dualepis)$
  induces a natural equivalence
  \begin{equation}
    \waddiagramsof{\partialmapsop[\factM]{\factcat}}
    \xlra{\simeq}
    \waddiagramsof{\factE^\op}
  \end{equation}
  for each weakly idempotent complete \pre{additive} \inftycat{} $\exwaddcat$.
\end{Cor}

The prototypical example of \autoref{cor:main-thm-for-partial-maps}
comes from the category $\finset$ of finite sets,
equipped with its canonical surjective-injective factorization system
$(\finsurj,\fininj)$;
in this case
$\partialmapsop{\finset}$ is precisely Segal's category
$\Gamma=\finsetpop$~\cite{Segal1974},
hence we get a natural equivalence
\begin{equation}
  \waddiagramsof{\Gamma}
  \xlra{\simeq}
  \waddiagramsof{\finsurj^\op}
\end{equation}
or, after dualizing (see \autoref{rem:main-thm-dualizes}),
\begin{equation}
  \label{eq:homotopy-Pirashvili}
  \waddiagramsof{\Gamma^\op}
  \xlra{\simeq}
  \waddiagramsof{\finsurj},
\end{equation}
for all weakly idempotent complete \pre{additive} \inftycats{} $\exwaddcat$.\\
We refer the reader to  \cite[Section~7]{LackStreet2015}
for many more examples in this spirit.

\begin{Rem}
  \label{rem:comparison-Helmstutler}
  Up to minor differences\footnote{
    For instance,
    Helmstutler's $\factM$ is not required
    to contain all isomorphisms.
    Unlike \nameLS{} (hence us) he does however require
    the pullback of an arrow in
    $\factE$
    along an arrow in
    $\factM$
    to lie in $\factE$ again;
    this amounts to saying that the set $M$ of \wmonos{}
    is closed under composition.
  },
  the pairs
  $(\partialmaps[\factM]{\factcat}, \factE)$
  arising from a factorization system
  $(\factE,\factM)$ as above
  are the \buzzword{conjugate pairs}
  $(\helcal{B},\helcal{A})$
  of Helmstutler\cite{Helmstutler2014}.
  For the convenience of the reader we provide
  a table translating Helmstutler's notation to the one of \nameLS{}
  (which we are using in this section):
  \begin{center}
    \begin{tabular}{c|c|c|c|c}
      \nameLS
      & $\factcat$
      & $\factM$
      & $\factE$
      & $\partialmaps[\factM]{\factcat}$
      \\
      \hline
      Helmstutler
      & $\helcal{U}$
      & $\helcal{I}$
      & $\helcal{A}$
      & $\helcal{B}$
    \end{tabular}
  \end{center}

  Helmstutler calls the arrows in
  $\factcat\subset\partialmaps[\factM]{\factcat}$
  \buzzword{regular}
  (because they are totally, and not just $\factM$-partially, defined)
  and the other arrows in
  $\partialmaps[\factM]{\factcat}$
  \buzzword{singular};
  this matches our use of those words.
  Moreover, he constructs a bimodule
  $
  \helcal{U}_+\colon \helcal{A}^\op\times\helcal{\localB}
  \lra\SetP
  $
  (which is precisely our bimodule
  $\regularsP$ from \autoref{cstr:wrapcat})
  and proves~\cite[Theorem 6.2]{Helmstutler2014} that it induces,
  for each left proper stable model category $\helcal{C}$,
  a Quillen equivalence
  $
  [\helcal{B}^\op,\helcal{C}]\adjarrows[\helcal{A}^\op,\helcal{C}]
  $
  (left adjoint on top).
  \autoref{cor:main-thm-for-partial-maps}
  is our version of this result,
  where we replace a Quillen equivalence of model categories
  by an equivalence of \inftycats{}.
  Note that the \self{duality} inherent to our \inftycategorical{} approach
  (see \autoref{rem:main-thm-dualizes})
  fixes the asymmetry problem addressed by Helmstutler
  in the note at the end of Section~6 in \cite{Helmstutler2014}.
\end{Rem}

\section{The proof}

\subsection{Cofinality lemmas}

\newcommand{\localXX}{X}
\newcommand{\localYY}{Y}

In order to get a better understanding of
the Kan extensions appearing in \autoref{thm:main-thm}
we use cofinality arguments to
simplify the relevant pointwise formulas.

\begin{Cstr}
  Fix an element $\locala\in\localA$.
  Consider the category
  \begin{equation}
    \localXX_{\locala}
    \coloneqq
    \dualepisof{\localaB}
    \pushout{\dualepisneqof{\localaB}}
    \rcone{\dualepisneqof{\localaB}},
  \end{equation}
  equipped with the functor
  $\localXX_{\locala}\to\overcat[\wrapcat]{\localBp}{\locala}$
  given by sending each
  $\localb\in\dualepisof{\localaB}$
  to the composition $\localb\ardualepi\localaB\to\locala$
  (which is the zero map for all $\localb\in\dualepisneqof{\localaB}$)
  and the cone point of
  $\rcone{\dualepisneqof{\localaB}}$
  to $\zerobj\to\locala$.
  Since $\dualepisneqof{\localaB}$ is
  (equivalent to) a poset,
  the same is true for $\localXX_{\locala}$;
  the latter poset arises from $\dualepisof{\localaB}$
  by adding one new element which is
  bigger than all elements of $\dualepisof{\localaB}$
  except its terminal object
  $\Id\colon\localaB\ardualepi\localaB$.\\
  Fix an element $\localb\in\localB$.
  Denote by
  $\localYY_{\localb}\subset \undercat[\wrapcat]{\localb}{(\localAp)}$
  the (\non{full}) subcategory defined as follows:
  \begin{itemize}
  \item
    objects are the maps
    $\localb\arepi \passtoA{\localb[1]}$
    corresponding to \wepis{} $\localb\arepi\localb[1]$ in $\localB$
    (recall that $\passtoA{\localb[1]}\in\localA$
    denotes the object corresponding to $\localb[1]\in\localB$)
    and the unique map $\localb\to\zerobj$.
  \item
    the only morphisms are isomorphisms under $\localb$
    and the zero morphisms $\passtoA{\localb[1]}\to \zerobj$.
  \end{itemize}
  Observe that $\localYY_{\localb}$ is equivalent to the right cone
  \begin{equation}
    \label{eq:finite-set-YY-b}
    \rcone{\set{\localb\arepi\localb[1]}}
  \end{equation}
  on the discrete set $\set{\localb\arepi\localb[1]}$
  containing some choice of representatives for the
  isomorphism classes of \wepis{} out of $\localb$;
  the cone point corresponds to the object
  $\localb\to\zerobj$ of $\localYY_{\localb}$.
\end{Cstr}

\begin{Lem}
  \label{lem:factorizations-XX}
  For each $\locala\in\localA$, the inclusion
  $\localXX_{\locala}\lhra\overcat[\wrapcat]{\localBp}{\locala}$
  is \htpyterminal{}.
\end{Lem}

Before we go into the rather technical proof of \autoref{lem:factorizations-XX},
we state the direct following corollary
which is what we will use going forward.

\begin{dCor}
  \label{cor:char-LKE-via-total-cofibers}
  \newcommand{\localLKE}{\exptdiagramB^1}
  Let $\exptdcat$ be a pointed \inftycat{} and
  $\exptdiagramB\colon\localBp\to\exptdcat$ a pointed diagram.
  Then a left Kan extension $\localLKE$ of $\exptdiagramB$
  along the inclusion $\localBp\hra\wrapcat$
  exists if and only if for each $\locala\in\localA$
  the total cofiber of the diagram
  \begin{equation}
    \dualepisof{\localaB}\lra \localBp\xra{\exptdiagramB} \exwaddcat
  \end{equation}
  exists in $\exwaddcat$.
  If it exists, this left Kan extension $\localLKE$
  is characterized pointwise at $\locala\in\localA$
  by the fact
  that it extends the diagram
  \begin{equation}
    \localXX_a
    \lra
    \overcat[\wrapcat]{\localBp}{\locala}
    \lra
    \localBp
    \xra{\exptdiagramB}
    \exwaddcat
  \end{equation}
  to a colimit cone in $\exwaddcat$
  with colimit $\localLKE(\locala)$.
\end{dCor}

\begin{Rem}
  \label{rem:char-LKE-if-the-colimit-exists}
  \newcommand{\localLKE}{\exptdiagramB^1}
  If the colimit of the diagram
  \begin{equation}
    \dualepisneqof{\localaB}\lra\localBp\xra{\exptdiagramB}\exwaddcat
  \end{equation}
  exists for each $\locala\in\localA$
  then we can characterize the left Kan extension
  as in \autoref{cor:char-LKE-via-total-cofibers}
  by the fact that it induces cofiber sequences
  \begin{equation}
    \label{eq:total-cofiber-diagram-a}
    \cdsquareNA[cC]
    {
      \colim\limits_{\localb\in\dualepisneqof{\localaB}}
      \localevalB\exptdiagramB{\localb}
    }
    {\localevalB\exptdiagramB\zerobj\simeq\zerobj}
    {\localevalB\exptdiagramB{\localaB}}
    {\localevalA\localLKE{\locala}}
  \end{equation}
  We will show in \autoref{prop:two-sided-Kan-extension}
  that this colimit always exists when
  $\exwaddcat$ is \waddjectives{}.
\end{Rem}

\begin{Prf}[of \autoref{lem:factorizations-XX}]
  \newcommand{\localf}{f}
  \newcommand{\locals}{s}
  \newcommand{\localx}{x}
  \newcommand{\localcatXfact}{
    \undercat[{\overcat[\wrapcat]{\localBp}{\locala}}]{\localu}{\localXX}
  }
  Fix $\locala\in\localA$ and let us abbreviate
  $\localXX\coloneqq\localXX_{\locala}$
  to avoid proliferating subscripts.
  For each $\localb\in\localB$ and each arrow
  $\localu\colon\localb\to\locala$ in $\wrapcat$,
  the undercategory $\localcatXfact$
  can be described explicitly as follows:
  \begin{itemize}
  \item
    objects are factorizations
    \begin{equation}
      \label{eq:factorization-XX}
      \cdtriangle
      {\localb}{\localb[1]}{\locala}
      {\localf}{\localx}{\localu}
    \end{equation}
    of $\localu$ in $\wrapcat$,
    where the arrow $\localx\colon\localb[1]\to\locala$
    is required to lie in $\localXX$;
  \item
    a morphism $(\localx,\localf)\to(\localx',\localf')$
    between such factorizations is simply an arrow
    $\localx\to\localx'$ in $\localXX$
    compatible with $\localf$ and $\localf'$.
  \end{itemize}
  Observe that $\localcatXfact$ is a poset (because $\localXX$ is).
  To prove that
  $\localXX\hra\overcat[\wrapcat]{\localBp}{\locala}$
  is \htpyinitial{} we have to show that all these categories of factorizations
  are \wcontractible{}.
  We distinguish two cases:
  \begin{itemize}
  \item
    Assume that $\localu\colon\localb\to\locala$ is a \non{zero}.
    Then the only factorization of $\localu$ through an object of $\localXX$
    is the tautological factorization
    $\localu\colon\localb\xra{\localuB}\localaB\xra{\canmapBA}\locala$
    hence the category $\localcatXfact$ is a singleton.
  \item
    Assume that $\localu\colon\localb\xra{\zeromap}\localu$ is the zero map.
    In this case there are three types of factorizations:
    \begin{enumerate}
      \item
        \label{item:XX-factorizations-1}
        given a \non{invertible} \wdualepi{}
        $\passtoB{\localx}\colon\localb[1]\ardualepineq\localaB$
        and given \emph{any} map
        $\localf\colon\localb\to\localb[1]$ in $\localBp$,
        there is a factorization
        $\zeromap\colon\localb\xra{\localf}\localb[1]\xra{\localx}\locala$;
      \item
        \label{item:XX-factorizations-2}
        for each \emph{singular} map
        $\locals\colon \localb\to\localaB$ in $\localBp$,
        there is a factorization
        $\zeromap\colon\localb\xra{\locals}\localaB\xra{\canmapBA}\locala$;
      \item
        \label{item:XX-factorizations-3}
        there is the zero factorization
        $\zeromap\colon\localb\xra{}\zerobj\xra{}\locala$.
      \end{enumerate}
      Denote by $(\localx,\localf)$, $(\canmapBA,\locals)$ and $0$
      the objects of $\localcatXfact$
      corresponding to the factorizations of type
      \ref{item:XX-factorizations-1},
      \ref{item:XX-factorizations-2} and
      \ref{item:XX-factorizations-3}, respectively.
      \newcommand{\XXubases}{Z_\locals}
      \newcommand{\XXubase}{Z}
      Denote by $\XXubase\subset\localcatXfact$
      the \sub{poset} consisting of the objects $(\localx,\localf)$
      of the first.
      For each singular map
      $\locals\colon \localb\to\localaB$,
      denote by $\XXubases\subset\XXubase$
      the \sub{poset} consisting of those $(\localx,\localf)$,
      where the composite
      $\localb\xra{\localf}\localb[1]\xra{\passtoB{\localx}}\localaB$
      is equal to $\locals\colon\localb\to\localaB$ .
      We now describe the morphisms in the category
      $\localcatXfact$.
      \begin{itemize}
      \item
        For each factorization $(\localx,\localf)$ with
        $\locals\coloneqq \passtoB{\localx}\circ\localf$
        (as maps $\localb\to\localaB$),
        we have a unique map $(\localx,\localf)\to(\canmapBA,\locals)$.
        There are no other maps between factorizations of types
        \ref{item:XX-factorizations-1} and \ref{item:XX-factorizations-2}.
        In other words,
        the \sub{poset}
        $\XXubases\cup\set{(\canmapBA,\locals)}\subset\localcatXfact$
        is a (right) cone on $\XXubases$ with maximum $(\canmapBA,\locals)$
      \item
        There are no maps between factorizations of types
        \ref{item:XX-factorizations-2} and \ref{item:XX-factorizations-3}
        (because there are no maps between
        $\localaB\to\locala$ and $\zerobj\to\locala$ in $\localXX$)
      \item
        For each factorization $(\localx,\localf)$,
        we have a unique map $(\localx,\localf)\to\zerobj$.
        There are no other maps between factorizations of types
        \ref{item:XX-factorizations-1} and \ref{item:XX-factorizations-3}.
        In other words,
        the \sub{poset}
        $\XXubase\cup\set{\zerobj}\subset\localcatXfact$
        is a (right) cone on $\XXubase$ with maximum $\zerobj$.
      \end{itemize}
      It follows that we have the following pushout of simplicial sets:
      \begin{equation}
        \label{eq:inprf:XX-attach-cones}
        \cdsquareOpt
        {
          \coprod\limits_{\locals\in\singulars(\localb,\localaB)}
          \nerve{\XXubases}
        }
        {
          \coprod\limits_{\locals\in\singulars(\localb,\localaB)}
          \rcone{\nerve{\XXubases}}
        }
        {\rcone{\nerve{\XXubase}}}
        {\nerve{\localcatXfact}}
        {hookrightarrow}{}{}{hookrightarrow}
      \end{equation}
      which is a (Kan) homotopy pushout
      because the top horizontal map is a monomorphism.
      By \ref{item:cooltech-fact-property},
      each singular arrow
      $\locals\colon\localb\to\localaB$
      admits a unique (up to unique isomorphism) factorization
      $\locals\colon\localb\to\localb[1]\ardualepineq\localaB$,
      where $\localb\to\localb[1]$ is regular
      and $\localb[1]\ardualepineq\localaB$
      is a \non{invertible} \wdualepi{};
      viewed as a factorization of
      $\zeromap\colon\localb\to\locala$
      it is an initial object of the category $\XXubases$,
      which is hence contractible.
      Therefore the top horizontal map
      in the square~\eqref{eq:inprf:XX-attach-cones}
      is a (Kan) weak equivalence,
      hence also the bottom horizontal map;
      this concludes the proof because
      $\rcone{\nerve{\XXubase}}$ is contractible.
      \qedhere
    \end{itemize}
\end{Prf}

\begin{Lem}
  \label{lem:RKE-by-products-initial}
  For each $\localb\in\localB$, the inclusion
  $\localYY_{\localb}\lhra \undercat[\wrapcat]{\localb}{(\localAp)}$
  is \htpyinitial{}.
\end{Lem}

\begin{Prf}
  \newcommand{\localf}{f}
  \newcommand{\locals}{s}
  \newcommand{\localy}{y}
  \newcommand{\localcatYfact}{
    \overcat{\localYY}{\localu}
  }
  Fix $\localb\in\localB$ and abbreviate $\localYY\coloneqq\localYY_{\localb}$.
  Similarly to the proof of \autoref{lem:factorizations-XX},
  we have to show that,
  for each $\locala\in\localAp$ and each map
  $\localu\colon\localb\to\locala$,
  the groupoid $\localcatYfact$ of factorizations
  \begin{equation}
    \label{eq:factorization-YY}
    \cdtriangle
    {\localb}{\localbA[1]}{\locala}
    {\localy}{\localf}{\localu}
  \end{equation}
  (with $\localy\in\localYY$) is \wcontractible{}.
  Again, we distinguish two cases:
  \begin{itemize}
  \item
    If the arrow $\localu\colon\localb\to\locala$ is \non{zero}
    then factorizations \eqref{eq:factorization-YY}
    are the same as $(\mcapr)\circ\epis$-factorizations
    of the corresponding regular map
    $\localuB\colon\localb\to\localaB$.
    By \ref{item:cooltech-fact-property},
    the groupoid of such factorizations is (equivalent to) a point.
  \item
    If the arrow $\localu$ is the zero
    then the zero factorization
    $\localu\colon\localb\to\zerobj\to\locala$
    is a terminal object of the category $\localcatYfact$,
    which is hence contractible.
    \qedhere
  \end{itemize}
\end{Prf}

\begin{Cor}
  \label{cor:RKE-by-products}
  \newcommand{\localRKE}{\exptdiagram}
  Let $\exwaddcat$ be a pointed \inftycat{} with finite products.
  Every pointed diagram
  $\exptdiagramA\colon\localAp\to\exwaddcat$
  admits a right Kan extension
  $\localRKE\colon\wrapcat\to\exwaddcat$
  along the inclusion
  $\localAp\hra\wrapcat$.
  Moreover,
  this right Kan extension is characterized pointwise by the product cones
  \begin{equation}
    \label{eq:product-formula-for-RKE}
    \localevalB\localRKE{\localb}\xra{\simeq}
    \prod_{\localb\arepi\locala}\localevalA\exptdiagramA{\locala}
  \end{equation}
  indexed by equivalence classes of \wepis{} out of $\localb$.
\end{Cor}

\begin{Prf}
  \newcommand{\localRKE}{\exptdiagram}
  Let $\exptdiagramA\colon\localAp\to\exwaddcat$ be a pointed diagram
  and fix an objects $\localb\in\localB$.
  By \autoref{lem:RKE-by-products-initial}
  we can compute the pointwise right Kan extension $\localRKE$
  of $\exptdiagramA$ along $\localAp\hra\wrapcat$
  at $\localb$ as the limit
  \begin{equation*}
    \localevalB\localRKE{\localb}
    \xra{\simeq}
    \lim
    {
      \left(
        \undercat[\wrapcat]{\localb}{(\localAp)}
        \to
        \localAp
        \xra{\exptdiagramA}
        \exwaddcat
      \right)
    }
    \xra{\simeq}
    \lim
    {
      \left(
        \localYY_{\localb}
        \to
        \localAp
        \xra{\exptdiagramA}
        \exwaddcat
      \right)
    }
    \simeq
    \lim
    {
      \left(
        \rcone{\set{\localb\arepi\locala}}
        \to
        \localAp
        \xra{\exptdiagramA}
        \exwaddcat
      \right)
    }.
  \end{equation*}
  This limit formula is the same as the product formula
  \eqref{eq:product-formula-for-RKE}
  because the value of
  $\exptdiagramA$ on the cone point of
  $\rcone{\set{\localb\arepi\locala}}$
  is $\localevalA\exptdiagramA\zerobj\simeq\zerobj$.
\end{Prf}

\subsection{Inductive construction in the reduced case}

\newcommand{\wrapcatll}[1]{\wrapcat_{<#1}}
\newcommand{\wrapcatleq}[1]{\wrapcat_{\leq#1}}
\newcommand{\dualwrapcatll}[1]{\dualwrapcat_{<#1}}
\newcommand{\dualwrapcatleq}[1]{\dualwrapcat_{\leq#1}}

Throughout this section we assume that the \cooltech{}
$\excooltech=(\localB,\epis,\dualepis)$
is reduced,
\ie that $\localB=\dualepis\circ\epis$
and hence $\localAp=\quotbyneq\localB$.
By applying
\autoref{cstr:wrapcat}
to the reduced \cooltech{} $\excooltech$
and to its dual $\excooltech^\op$,
we obtain two categories
\begin{equation}
  \wrapcat
  =\wrapcat(\excooltech)
  \coloneqq
  \left(
    \begin{matrix}
      \localApv{1} & \regularsP\\
      \zerobj & \localBp
    \end{matrix}
  \right)
  \intxt{and}
  \dualwrapcat
  \coloneqq
  {\wrapcat(\excooltech^\op)}^\op
  =
  \left(
    \begin{matrix}
      \localBp & \coregularsP\\
      \zerobj & \localApv{0}
    \end{matrix}
  \right)
\end{equation}
where $\localApv{0}$ and $\localApv{1}$
are both just (a copy of) $\localAp$,
decorated with superscripts $0$ and $1$ to avoid confusing them.
For every $\locala\in\localA$
we denote by $\localav{0}$ its copy in
$\localAv{0}\subset{\dualwrapcat}$
and by $\localav{1}$ its copy in
$\localAv{1}\subset{\wrapcat}$.
Furthemore, we denote by $\wrapcatleq{\locala}\subset\wrapcat$
the full \sub{categories} spanned by $\localBp$
and by all the objects $\localav[1]{1}$ with
$\locala[1]\leq\locala$;
similarly, $\dualwrapcatleq{\locala}\subset\dualwrapcat$
is the full \sub{category} which contains $\localBp$
and all the objects $\localav[1]{0}$ with
$\locala[1]\leq\locala$.

\begin{Prop}
  \label{prop:two-sided-Kan-extension}
  Let $\exwaddcat$ be a \waddjectives{} \inftycat{} $\exwaddcat$
  and let $\exptdiagramB\colon\localBp\to\exwaddcat$
  be a pointed functor.
  Then there exist functors
  \begin{equation}
    \localXtension{0}\colon\dualwrapcat\to\exwaddcat
    \intxt{and}
    \localXtension{1}\colon\wrapcat\to\exwaddcat
  \end{equation}
  which are right and left Kan extension of $\exptdiagramB$, respectively.
  Moreover the functors $\localXtension{0}$ and $\localXtension{1}$
  are a left Kan extension and a right Kan extension
  of their restriction to $\localApv{0}$ and $\localApv{1}$,
  respectively.
\end{Prop}

\begin{Rem}
  By \autoref{cor:RKE-by-products},
  the \qquote{moreover} part of
  \autoref{prop:two-sided-Kan-extension}
  is saying that
  for each $\localb\in\localB$
  the diagrams $\localXtension{0}$ and $\localXtension{1}$
  induce direct sum decompositions
  \begin{equation}
    \label{eq:prod-coprod-decomp-of-b}
    \coprod_{\locala\ardualepi\localb}\localXtension{0}(\locala)
    \xra{\simeq}
    \localXtension{0}(\localb)
    =
    \localevalB\exptdiagramB{\localb}
    \intxt{and}
    \localevalB\exptdiagramB{\localb}
    =
    \localXtension{1}(\localb)
    \xra{\simeq}
    \prod_{\localb\arepi\locala}\localXtension{1}(\locala)
  \end{equation}
    where the coproduct/product is indexed over equivalence classes
    of \wdualepis{} into $\localb$ and
    \wepis{} out of $\localb$, respectively.
\end{Rem}

\begin{Rem}
  \label{rem:decompositions-are-adapted}
  It follows from the universal property of the coproduct
  that each \wdualepi{}
  $\localb[1]\ardualepi\localb$
  induces a commutative square
  \begin{equation}
    \cdsquareOpt
    {\coprod\limits_{\locala\ardualepi\localb[1]}\localXtension{0}(\locala)}
    {\localevalB\exptdiagramB{\localb[1]}}
    {\coprod\limits_{\locala\ardualepi\localb}\localXtension{0}(\locala)}
    {\localevalB\exptdiagramB{\localb}}
    {"\simeq"}
    {}
    {}
    {"\simeq"}
  \end{equation}
  where the left vertical map is the inclusion of those summands
  that are labeled by a \wdualepi{} which factors through
  $\localb[1]\ardualepi\localb$.
  Similarly each \wepi{} $\localb\arepi\localb[1]$
  induces projection onto those factors of the decomposition
  $\localevalB\exptdiagramB{\localb}\simeq
  \prod_{\localb\arepi\locala}\localXtension{1}(\locala)$
  that are indexed by \wepis{} which factor through
  $\localb\arepi\localb[1]$.
\end{Rem}

\begin{Prf}
  \newcommand{\localXtensionll}[2]{\exptdiagramB^{#1}_{<#2}}
  \newcommand{\localXtensionleq}[2]{\exptdiagramB^{#1}_{\leq#2}}
  \newcommand{\localXhbig}{\ol X}
  \newcommand\localcomparisonmap{\Phi}
  \newcommand{\localmatrixentry}{\Phi_{\locala[2],\locala[1]}}
  \newcommand{\localmatrixcomponents}{
    \matrixcomponents[\localaB]{\localaB[2]}{\localaB[1]}
  }
  For each $\locala\in\localA$ we prove:
  \begin{enumerate}
  \item
    \label{item:RKE-Vdualleqa}
    A right Kan extension $\localXtensionleq{0}{\locala}$
    of $\exptdiagramB$ along $\localBp\hra\dualwrapcatleq{\locala}$ exists.
  \item
    \label{item:LKE-Vleqa}
    A left Kan extension $\localXtensionleq{1}{\locala}$
    of $\exptdiagramB$ along $\localBp\hra\wrapcatleq{\locala}$ exists.
  \item
    \label{item:Wleqa-direct-sums}
    Each choice of such Kan extensions
    $\localXtensionleq{0}{\locala}$ and $\localXtensionleq{1}{\locala}$
    induces, for each $\localb\leq\localaB$,
    direct sum decompositions as in \eqref{eq:prod-coprod-decomp-of-b};
    moreover, the composition
    \begin{equation}
      \localXtensionleq{0}{\locala}(\localav{0})
      \lra
      \localevalB\exptdiagramB{\localaB}
      \lra
      \localXtensionleq{1}{\locala}(\localav{1})
    \end{equation}
    is an equivalence in $\exwaddcat$.
  \end{enumerate}
  By induction on the number
  $
  \card{\pi_0\dualepisof{\locala}}
  =
  \card{\pi_0\episof{\locala}}
  $
  we may assume that we have proved
  \ref{item:RKE-Vdualleqa},
  \ref{item:LKE-Vleqa} and
  \ref{item:Wleqa-direct-sums}
  for all objects of $\localA$ which are strictly smaller than $\locala$.
  Choose a right Kan extension
  $\localXtensionll{0}{\locala}\colon\dualwrapcatll{\locala}\to\exwaddcat$
  and a left Kan extension
  $\localXtensionll{1}{\locala}\colon\wrapcatll{\locala}\to\exwaddcat$
  of $\exptdiagramB\colon\localBp\to\exwaddcat$
  (they exist pointwise by assumption).
  By assumption, $\localXtensionll{0}{\locala}$
  induces coproduct decompositions
  $\coprod\limits_{\locala[1]\ardualepi\localb}
  \localXtensionll{0}{\locala}(\locala[1])
  \xra{\simeq}
  \localevalB\exptdiagramB{\localb}$
  for all $\localb<\localaB$.
  Since all \wdualepis{} induce compatible inclusions of summands
  (see \autoref{rem:decompositions-are-adapted}),
  the diagram $\localXtension{0}$ provides an identification
  \begin{equation}
   \label{eq:express-colimit-via-coproduct}
    \coprod_{\locala[1]\ardualepineq\localaB}
    \localXtensionll{0}{\locala}(\locala[1])
    \simeq
    \colim_{\localb\in\dualepisneqof{\localaB}}\localevalB\exptdiagramB{\localb}
  \end{equation}
  where the coproduct is indexed over
  equivalence classes of \emph{\non{invertible}} \wdualepis{};
  moreover, this identification \eqref{eq:express-colimit-via-coproduct}
  is compatible with the respective structure maps to
  $\localevalB\exptdiagramB{\localaB}$.
  By applying the dual argument to
  $\localXtensionll{1}{\locala}\colon\wrapcatll{\locala}\to\exwaddcat$
  we obtain an identification
  \begin{equation}
    \lim_{\localb\in\episneqof{\localaB}}\localevalB\exptdiagramB{\localb}
    \simeq
    \prod_{\localaB\arepineq\locala[1]}
    \localXtensionll{1}{\locala}(\locala[1]),
  \end{equation}
  again compatible with the structure maps from
  $\localevalB\exptdiagramB{\localaB}$.
  We analyze the two composable maps
  \begin{equation}
    \label{eq:candidate-s-r-coprod-prod}
    \coprod_{\locala[1]\ardualepineq\localaB}
    \localXtensionll{0}{\locala}(\locala[1])
    \lra
    \localevalB\exptdiagramB\localaB
    \lra
    \prod_{\localaB\arepineq\locala[2]}
    \localXtensionll{1}{\locala}(\locala[2])
  \end{equation}
  and their composite $\localcomparisonmap$
  in terms of the components
  $\localmatrixentry
  \colon
  \localXtensionll{0}{\locala}(\locala[1])
  \to
  \localevalB\exptdiagramB\localaB
  \to
  \localXtensionll{1}{\locala}(\locala[2])
  $.
  We have the commutative diagram in
  $\dualwrapcatll{\locala}\pushout{\localBp}\wrapcatll{\locala}$
  \begin{equation}
    \label{eq:butterfly-aaa}
    \begin{tikzcd}
      & \localaB[2]\ar[rd,"\canmapBA"]&\\
      \localav[1]{0}\ar[r]\ar[rd,"\canmapAB"']
      & \localaB\ar[r]\ar[u,twoheadrightarrow]
      & \localav[2]{1}\\
      & \localaB[1]\ar[u,rightarrowtail]&
    \end{tikzcd}
  \end{equation}
  where the vertical morphisms are the \wdualepi{}
  $\localaB[1]\ardualepi\localaB$ and the \wepi{}
  $\localaB\arepi\localaB[1]$;
  their composition is---by definition---the map
  ${\localmatrixcomponents}.$
  Therefore, the map $\localmatrixentry$ is equivalent to the composition
  \begin{equation}
    \localmatrixentry
    \colon
    \localXtensionll{0}{\locala}(\localav[1]{0})
    \lra
    \localevalB\exptdiagramB{\localaB[1]}
    \xra{
      \exptdiagramB({\localmatrixcomponents})
    }
    \localevalB\exptdiagramB{\localaB[2]}
    \lra
    \localXtensionll{1}{\locala}(\localav[2]{1}).
  \end{equation}
  It follows that:
  \begin{itemize}
  \item
    If $\localmatrixcomponents$
    is an isomorphism in $\localBp$
    (without loss of generality, the identity)
    then $\localmatrixentry$ is an equivalence
    by the induction hypothesis \ref{item:Wleqa-direct-sums};
  \item
    If $\localmatrixcomponents$
    is not an isomorphism in $\localBp$
    then it must be either singular or cosingular.
    If it is singular then the composition
    $\localaB[1]\to\localaB[2]\to\localav[2]{1}$
    factors through $\zerobj\in\wrapcat$;
    if it is cosingular then the composition
    $\localav[1]{0}\to\localaB[1]\to\localaB[2]$
    factors through $\zerobj\in\dualwrapcat$;
    in either case $\localmatrixentry$ factors through
    $\localevalB\exptdiagramB{\zerobj}\simeq\zerobj$.
  \end{itemize}
  Therefore it follows from \ref{item:cooltech-UT-matrix}
  that $(\localmatrixentry)$
  is an upper triangular matrix with invertible diagonal entries;
  hence $\localcomparisonmap$ is invertible because $\exwaddcat$ is additive
  (see \autoref{lem:upper-triangles-are-invertible}).
  This means that the two composable maps \eqref{eq:candidate-s-r-coprod-prod}
  are a section-retraction pair.
  Since $\exwaddcat$ is weakly idempotent complete,
  this section-retraction pair admits a complement,
  \ie there is an essentially unique diagram
  \newcommand{\localK}{K}
  \newcommand{\localQ}{Q}
  \begin{equation}
   \label{eq:add-kernels-cokernels-to-section-retraction}
   \cdcomplsecrecNA
   {\localevalB\exptdiagramB\localaB}
   {
     \coprod\limits_{\locala[1]\ardualepineq\localaB}
     \localXtensionll{0}{\locala}(\locala[1])
   }
   {
     \prod\limits_{\localaB\arepineq\locala[2]}
     \localXtensionll{1}{\locala}(\locala[2])
   }
   {\localK}
   {\localQ}
  \end{equation}
  where all squares are biCartesian.
  By
  \ref{cor:char-LKE-via-total-cofibers}
  (or, more precisely, by \ref{rem:char-LKE-if-the-colimit-exists})
  and the identification
  \eqref{eq:express-colimit-via-coproduct},
  we conclude that the pointwise left Kan extension
  $\localXtension{1}(\locala)$ of $\exptdiagramB$ at $\localav{1}$ exists
  and that its value on the structure map
  $\canmapBA\colon\localaB\to\localav{1}$
  is equivalent to the projection
  $\localevalB\exptdiagramB{\localaB}\to\localQ$.
  By the dual argument,
  we conclude that the pointwise right Kan extension
  $\localXtension{0}(\locala)$
  of $\exptdiagramB$ at $\localav{0}$ exists
  and that its value on the structure map
  $\canmapAB\colon\localav{0}\to\localaB$
  is equivalent to the inclusion
  $\localK\to\localevalB\exptdiagramB{\localaB}$.
  To establish the inductive step for property \ref{item:Wleqa-direct-sums},
  note that the
  diagram~\eqref{eq:add-kernels-cokernels-to-section-retraction}
  encodes the required coproduct decompositions
  \begin{equation}
    \coprod\limits_{\locala[1]\ardualepi\localaB}
    \localXtension{0}(\locala[1])
    =
    \localK\sqcup
    \coprod\limits_{\locala[1]\ardualepineq\localaB}
    \localXtensionll{0}{\locala}(\locala[1])
    \xra{\simeq}
    \localevalB\exptdiagramB\localaB
  \end{equation}
  (and similarly the required product decomposition)
  and the fact that the composition
  \begin{equation}
    \localK
    \simeq
    \localXtension{0}(\locala)
    \xra{\localXtension{0}(\canmapAB)}
    \localevalB\exptdiagramB\localaB
    \xra{\localXtension{1}(\canmapBA)}
    \localXtension{1}(\locala)
    \simeq
    \localQ
  \end{equation}
  is an equivalence.
\end{Prf}

\begin{Rem}
  \label{rem:when-matrix-in-prop-is-diagonal}
  \newcommand{\localmatrixentry}{\Phi_{\locala[2],\locala[1]}}
  If the \cooltech{} $\excooltech$ is \really{}
  then the matrix $(\localmatrixentry)_{\locala[2],\locala[1]}$
  is actually a diagonal matrix.
  Hence to invert it, we do not need $\exwaddcat$ to be additive
  but only \pre{additive}.
\end{Rem}

\begin{Rem}
  \label{rem:encoded-complementary-sec-rec-pair}
  From the proof of \autoref{prop:two-sided-Kan-extension}
  we can extract more detailed information.
  For each $\locala\in\localA$,
  the extensions
  $\localXtension{0}$ and $\localXtension{1}$
  encode two complementary section-retraction pairs
  \begin{equation}
    \cdcomplsecrecNA
    {\localevalB{\exdiagramB}{\localaB}}
    {\localXtension{0}(\locala)}
    {\localXtension{1}(\locala)}
    {
      \colim\limits_{\localb\in\dualepisneqof{\localaB}}
      \localevalB{\exdiagramB}{\localb}
    }
    {
      \lim\limits_{\localb\in\episneqof{\localaB}}
      \localevalB{\exdiagramB}{\localb}
    }
  \end{equation}
  (in particular, the indicated limits/colimits exist).
\end{Rem}

\subsection{Proof of the main theorem}

\begin{Prf} [of \autoref{thm:main-thm}]
  \newcommand{\localtraf}{\alpha}
  We first prove
  \ref{it:main-Kan-extensions-exist},
  \ref{it:main-master-equivalence},
  and
  \ref{it:main-colim-lim-formula}
  in the case where
  the \cooltech{} $(\localB,\epis,\dualepis)$ is reduced.
  In this case, we have the following ingredients:
  \begin{itemize}
  \item
    \autoref{cor:RKE-by-products}
    guarantees that the right Kan extension functor
    $\RKE\colon\waddiagramsPof{\localAp}
    \to\waddiagramsPof{\wrapcat}$
    exists.
    Moreover, the explicit formula \eqref{eq:product-formula-for-RKE}
    implies that for any natural transformation
    $\localtraf\colon\exptdiagramA'\to\exptdiagramA$
    of pointed diagrams $\localAp\to\exwaddcat$,
    the component
    $\localevalA\localtraf{\locala}
    \colon
    \localevalA{\exptdiagramA'}{\locala}
    \to
    \localevalA{\exptdiagramA}{\locala}$
    at $\locala\in\localA$
    is a factor of the corresponding right Kan extended transformation
    (with the notation as in \ref{cor:RKE-by-products})
    \begin{equation}
      \cdsquare
      {\localevalB{\exptdiagramB'}{\localaB}}
      {\localevalB{\exptdiagramB}{\localaB}}
      {\prod\limits_{\localaB\arepi\locala[1]}\localevalA{\exptdiagramA'}{\locala[1]}}
      {\prod\limits_{\localaB\arepi\locala[1]}\localevalA\exptdiagramA{\locala[1]}}
      {\localevalB{\RKE(\localtraf)}{\localaB}}
      {\simeq}
      {\simeq}
      {\prod\localevalA\localtraf{\locala[1]}}
    \end{equation}
    at $\localaB\in\localB$;
    hence it follows from \ref{lem:factors-of-equis-are-equis}
    that the composition
    \begin{equation}
      \waddiagramsPof{\localAp}
      \xra{\RKE}
      \waddiagramsPof{\wrapcat}
      \xra{\Res}
      \waddiagramsPof{\localBp}
    \end{equation}
    is conservative, \ie reflects equivalences.
  \item
    \autoref{prop:two-sided-Kan-extension} states in particular---%
    if we focus only on the statements about $\wrapcat$
    and not about $\dualwrapcat$---%
    that
    \begin{itemize}
    \item
      the left Kan extension functor
      $\LKE\colon\waddiagramsPof{\localBp}
      \to\waddiagramsPof{\wrapcat}$
      exists and
    \item
      on the image of this functor $\LKE$,
      the unit
      $\Id_{\waddiagramsPof{\wrapcat}}\to\RKE\circ\Res$
      of the adjunction
      \begin{equation}
        \Res\colon
        \waddiagramsPof{\wrapcat}
        \adjarrows
        \waddiagramsPof{\localAp}
        \noloc\RKE
      \end{equation}
      is an equivalence.
    \end{itemize}
    Since left Kan extension along the fully faithful functor
    $\localBp\hra\wrapcat$ is fully faithful,
    the unit
    $\Id_{\waddiagramsPof{\localBp}}\to\Res\circ\LKE$
    is  an equivalence.
    We conclude that the unit
    \begin{equation}
      \Id_{\waddiagramsPof{\localBp}}
      \lra
      \Res\circ\LKE
      =
      \Res\circ\Id_{\waddiagramsPof{\wrapcat}}\circ\LKE
      \lra
      \Res\circ\RKE\circ\Res\circ\LKE
    \end{equation}
    of the composite adjunction \ref{eq:master-span-LRKE}
    is also an equivalence.
  \end{itemize}
  This already proves
  \ref{it:main-Kan-extensions-exist};
  assertion
  \ref{it:main-master-equivalence}
  follows from the general fact about adjunctions that
  if the right adjoint is conserative and the unit is an equivalence
  then the whole adjunction is an adjoint equivalence.
  Assertion
  \ref{it:main-colim-lim-formula}
  is spelled out in
  \autoref{rem:encoded-complementary-sec-rec-pair}
  since
  $\localevalA\exptdiagramA{\locala}$
  is by definition equivalent to
  $\localXtension{1}(\locala)$.

  To prove
  \ref{it:main-Kan-extensions-exist},
  \ref{it:main-master-equivalence} and
  \ref{it:main-colim-lim-formula}
  when $\excooltech$ is not necessarily reduced,
  we make the following key observation:
  \begin{itemize}
  \item
    the criterion for constructing and detecting left Kan extension along
    $\localBp\hra\wrapcat$
    (\autoref{cor:char-LKE-via-total-cofibers})
    and the criterion for constructing and detecting right Kan extension along
    $\localAp\hra\wrapcat$
    (\autoref{cor:RKE-by-products})
    both only depend on the values of
    a diagram on the \wdualepis{} $\dualepis$ and on the \wepis{} $\epis$.
  \end{itemize}
  Therefore we can reduce to the reduced case (no pun intended)
  by replacing the original \cooltech{}
  with the reduced \cooltech{}
  \begin{equation}
    \ol{\excooltech}
    \coloneqq
    ( \dualepis\circ\epis
    , \epis
    , \dualepis
    ).
  \end{equation}

  To prove
  \ref{it:main-naturality},
  note that the right Kan extension
  $\RKE\colon\waddiagramsPof{\localAp}\lra\waddiagramsPof{\wrapcat}$
  is natural in $\exwaddcat$ with respect to all functors which
  preserve the relevant pointwise limits;
  since all these pointwise limits are just products,
  this is true for every additive functor.
\end{Prf}

\section{Comparison with...}

\subsection{...the setting of \nameLS}

\label{sec:comparison-Lack-Street}

We provide a short dictionary/comparison between
our setup described in \autoref{sec:setup} and \autoref{sec:key-constructions}
and the setting of \nameLS{}~\cite[Section~2]{LackStreet2015}.
Unless stated otherwise, references in this section refer to
their revised arXiv paper\cite{LackStreet2014},
\emph{not} the published one~\cite{LackStreet2015}
(see also the corrigendum~\cite{LackStreet2015corr});
we freely use the notation of \cite[Section~2]{LackStreet2014}.

Their category $\LSscr{P}$ is the \emph{dual} of our category $\localB$.
Under this duality we have the following table of correspondence:

\begin{center}
  \begin{tabular}{c|c|c|c|c|c|c|c}
    \nameLS
    & $\LSscr{P}$
    &  $\LSscr{M}$
    &  $\LSscr{M}^{\ast}$
    &  $\LSscr{R}$
    &  $\LSscr{D}$
    &  $\LSscr{S}$
    &  $\setP{u}{s_u\in\LSscr{R}}$
    \\
    \hline
    our setup
    &  $\localB$
    &  $\epis$
    &  $\dualepis$
    &  $\mcapr$
    &  $\localA$
    &  $\monos$
    &  $\regulars$
  \end{tabular}
\end{center}

\nameLS{} take as part of the data an isomorphism
$(\blank)^\ast\colon\LSscr{M}^\op\cong \LSscr{M}^\ast$
(which in our language would be written as
$(\blank)^\dual\colon\epis^\op\cong\dualepis$)
which is the identity on objects and
satisfies $m^\ast\circ m=\Id$ for all arrows $m$ in $\LSscr{M}$.
Their Assumption 2.5 translates to the fact that the set
$\pi_0\episof{\localb}$ is finite for each $\localb\in\localB$;
Assumption 2.6 is saying
that for each $\localb\in\localB$
there exists a linear order on $\pi_0\episof{\localb}$
such that the matrix
$
\matrixcomponents[\localb]{}{(\blank)^\dual}\colon
\pi_0\episof{\localb}\times\pi_0\episof{\localb}
\ra
\pi_0\Ar{\localB}$
has only singular entries below the diagonal.
In our setup,
\ref{item:cooltech-UT-matrix} replaces all these ingredients and
repackages them as a property which more directly reflects the final use:
what we ultimately want to exploit is that
certain unipotent upper triangular matrices
\eqref{eq:candidate-s-r-coprod-prod}
induced from the matrices $\matrixcomponents[\localb]{}{}$
can be inverted in any additive \inftycat{}.
Note that while \nameLS{} require
the matrix entries below the diagonal to be singular,
it suffices for our purposes if they are \non{invertible}.\\
Furthermore:
\begin{itemize}
\item
  Their Assumption 2.1 and Assumption 2.4
  correspond precisely to our axioms
  \ref{item:cooltech-fact-property}
  and
  \ref{item:cooltech-epis-dual-epis-cat},
  respectively.
\item
  Their Assumption 2.2 translates to our axiom
  \ref{item:cooltech-monos-cat}.
\item
  Their Assumption 2.3 translates to
  $(\monos\cap\regulars)\circ
  \dualepisneq
  \subset \localB\setminus(\monos\cap\regulars)$
  and is, \latin{a priori}, weaker than our axiom
  \ref{item:cooltech-sing-ideal}.
  However, they use Assumption 2.3
  (in the presence of the other assumptions)
  to prove Proposition 2.10(b)
  which states that if two composable arrows
  $v,u$ satisfy $s_{v}\not\in\LSscr{R}$ and $u\in \LSscr{S}$,
  then also $s_{vu}\not\in\LSscr{R}$.
  This statement translates to
  $\monos\circ\singulars\subseteq\singulars$,
  which is precisely
  \ref{item:cooltech-sing-ideal}.
\end{itemize}
The preceding discussion proves:
\begin{dCor}
  Let $\LSscr{P}$, $\LSscr{M}$, $\LSscr{M}^\ast$ and $\LSscr{D}$
  be as in \cite[Section~2]{LackStreet2014}.
  Then
  $
  \excooltech=
  (\LSscr{P}^\op
  ,\LSscr{M}^\op
  ,{(\LSscr{M}^\ast)}^\op
  )
  $
 is a \cooltech{} with associated pointed category
  $\localApof{\excooltech}=\quotcat{\LSscr{D}^\op}{\zeroideal}$.
\end{dCor}

The main tool in the proof by \nameLS{}
is what they call the \buzzword{kernel module}
\cite[Section~4]{LackStreet2015}
\begin{equation}
  M\colon\LSscr{D}^\op\times\LSscr{P}\lra 1/\mathrm{Set}
\end{equation}
(where $1/\mathrm{Set}$ is their notation for the category of pointed sets);
it corresponds to our
$\localAp$-$\localBp$-bimodule
\begin{equation}
  \regularsP
  \colon
  \localBp^\op\times\localAp\lra\SetP
\end{equation}
which we encode in its upper triangular category $\wrapcat$.
Their main theorem
~\cite[Theorem~6.7]{LackStreet2014}%
~\cite[Theorem~6.8]{LackStreet2015}
states that for each idempotent complete additive $1$-category $\LSscr{X}$,
the kernel module $M$ induces an equivalence
\begin{equation}
  \Fun(\LSscr{P},\LSscr{X})
  \simeq
  \Fun_{\SetP}(\LSscr{D},\LSscr{X})
\end{equation}
where $\Fun_{\SetP}$ denotes the category $\SetP$-enriched functors.
Instead of using $\SetP$-enriched categories
(or rather $\SpacesP$-enriched \inftycats{};
see also \autoref{rem:what-about-pointed-enriched})
we chose to work with pointed categories
and phrase our main result in terms of pointed functors
on $\localAp=\quotcat{\localA}{\zeroideal}$.
Therefore \autoref{cor:main-thm-1-categorical} recovers their result because,
for each pointed $1$-category $\exptdordcat$
and each $\SetP$-enriched category $\localA$,
the inclusion
$\localA\hra\localAp$
induces an equivalence of categories
$
\ptdFun{\localAp}{\exptdordcat}
\xra{\simeq}
\Fun_{\SetP}(\localA,\exptdordcat)
$
(see \autoref{rem:1-trivializing-ff}).

\subsection{...Lurie's stable Dold--Kan correspondence}

\label{sec:comparison-Lurie}

Let $\exstablecat$ be an \inftycat{} with finite colimits
and consider the functor
\begin{equation}
  \label{eq:Lurie-DK-functor}
  \Fun(\Dop,\exstablecat)\lra\Fun(\BN,\exstablecat),
\end{equation}
which sends a simplicial object
$\exsobj\colon\Dop\to\exstablecat$
to the filtered object
\begin{equation}
  \exfilteredDK\colon
  \colim\exsobjleq{0}
  \lra\cdots\lra
  \colim\exsobjleq{\nn-1}
  \lra
  \colim\exsobjleq{\nn}
  \lra\cdots
\end{equation}
of its partial colimits
$\exfilteredDK_{\nn}\coloneqq\colim\exsobjleq{\nn}
=\colim(\exsobjleq{\nn}\colon\Dopleq{\nn}\hra\Dop\xra{\exsobj}\exstablecat)$.
Lurie's stable \DKcorr{}~\cite[Theorem 1.2.4.1]{Lurie2017}
states that the functor \eqref{eq:Lurie-DK-functor} is
an equivalence
when the target $\exstablecat$ is a \emph{stable} \inftycat{}.
The functor \eqref{eq:Lurie-DK-functor}
lifts the ordinary \DKcorr{} in the following sense:
Each filtered object $\exfilteredDK$
in a stable \inftycat{} $\exstablecat$ gives rise
to a connective chain complex
\begin{equation}
  \label{eq:stable-incoherent-chain-complex}
  \htpycat\extermcplx_0
  \lla\cdots\lla
  \htpycat\extermcplx_{\nn-1}
  \lla
  \htpycat\extermcplx_{\nn}
  \lla\cdots
\end{equation}
in the homotopy category $\htpycat{\exstablecat}$,
with
$\extermcplx_{\nn}\coloneqq
\loops^\nn\cof(\exfilteredDK_{\nn-1}\to{\exfilteredDK_{\nn}})$.
Moreover, there is a commutative diagram
\begin{equation}
  \label{eq:comparison-Lurie-DK}
  \begin{tikzcd}
    &{\Fun(\BN,\exstablecat)}\ar[dd, bend left=60]\ar[d,dotted]
    &\exfilteredDK\ar[d,mapsto,dotted]\ar[dd, bend left=60,mapsto]
    \\
    {\Fun(\Dop,\exstablecat)}
    \ar[ur,"\simeq"]\ar[r,"\simeq",leftrightarrow]\ar[d]
    &\connChof{\exstablecat}\ar[d]
    &\extermcplx\ar[d,mapsto]
    \\
    {\Fun(\Dop,\htpycat\exstablecat)}\ar[r,"\simeq",leftrightarrow]
    &{\connChof{\htpycat{\exstablecat}}}
    &\htpycat\extermcplx
  \end{tikzcd}
\end{equation}
where the top diagonal functor is
\eqref{eq:Lurie-DK-functor}
and the lower commutative square is the naturality square of
\autoref{rem:naturality-to-htpy-cat}.
In particular,
the dotted equivalence $\exfilteredDK\mapsto\extermcplx$
exists and functorially lifts the incoherent
chain complex~\eqref{eq:stable-incoherent-chain-complex}
to a coherent one.

If we only assume that the target $\exstablecat$
is weakly idempotent complete additive
but not necessarily stable
then,
even if sufficient colimits exist
to define the functor \eqref{eq:Lurie-DK-functor},
it need not be an equivalence anymore;
similarly, the dotted functor
$\exfilteredDK\mapsto\extermcplx$
(or even $\exfilteredDK\mapsto\htpycat{\extermcplx}$)
does not exist in this generality.
For instance,
in the \inftycat{} of \emph{connective} spectra
the filtered object
$0\to \spherespectrum\to \spherespectrum\to \spherespectrum\to\dots$
(which would correspond to the chain complex
$0\la \spherespectrum\shift{-1} \la0 \la 0\la$)
does not arise from a simplicial object.

\begin{Rem}
  A systematic study of the relationship between
  coherent chain complexes and
  filtered objects in stable \inftycats{}
  is part of \SAriotta{}'s \PhD thesis~\cite{Ariotta}.
  In particular, he directly constructs
  an equivalence $\Fun(\BN,\exstablecat)\simeq\connChof{\exstablecat}$
  of \inftycats{}
  which we expect to agree with
  the vertical dotted equivalence in \eqref{eq:comparison-Lurie-DK}
  obtained by combining our result with Lurie's.
\end{Rem}

\section{Application: measuring Kan extensions}

Let $\excooltech=(\localB,\epis,\dualepis)$ be a \cooltech{}
with associated quotient
$\localAp=\localApof{\excooltech}$.
Let
$\exdiagramB\colon\localB\to\exwaddcat$
be a diagram in a \waddjectives{} \inftycat{} $\exwaddcat$ and let
$\exptdiagramA\colon\localAp\to\exwaddcat$
be the pointed functor corresponding to $\exdiagramB$
under the equivalence of
\autoref{cor:main-thm-without-pt}.

In this section, we set out to answer the following question:

\begin{Qstn}
  \label{qstn:what-does-A-know-about-B}
  What do the values of the diagram
  $\exptdiagramA\colon\localAp\to\exwaddcat$
  tell us about the original diagram
  $\exdiagramB\colon\localB\to\exwaddcat$?
\end{Qstn}

The rough answer is that in favorable situations
$\exptdiagramA$ \roughly{measures}
how far away $\exdiagramB$ is from being a Kan extension of its restriction
$\exdiagramBll{\nn}$.
To make this precise,
we make the following definition:

\begin{Def}
  The \cooltech{} $\excooltech$ is called \introduce{\coolsmonotone{}}
  if all \wmonos{} make objects bigger,
  \ie if we have $\localb[1]\leq\localb$
  whenever there exists a \wmono{}
  $\exmono\colon \localb[1]\to\localb$.
  We say that $\excooltech$ is \introduce{\coolpmonotone{}}
  if \wmonos{} at least do not make objects smaller,
  \ie there are no \wmonos{} $\localb[1]\to\localb$
  if $\localb[1]>\localb$.
\end{Def}

\begin{Rem}
  If the partial order $\leq$ on $\pi_0\localB$ is total,
  then the notions of \coolsmonotone{} and \coolpmonotone{} agree;
  in general being \coolpmonotone{} is weaker than being \coolsmonotone{}.
\end{Rem}

\begin{Rem}
  If $(\pi_0\localB,\leq)$ is the
  poset $(\BN,\leq)$ of natural numbers,
  then each \coolsmonotone{} \cooltech{} is a generalized Reedy category
  in the sense of
  \cite[Definition~1.1]{BergerMoerdijk2011}
  or \cite[Definition~8.1.1]{Cisinski2006}
  with tautological degree function $\Ob\localB\to\pi_0\localB\cong\BN$,
  degree-raising arrows $\monos$ and degree-lowering maps $\epis$.
\end{Rem}

\begin{Rem}
  Whether the \cooltech{}
  $(\localB,\epis,\dualepis)$
  is \pcoolpmonotone{} does not depend on $\dualepis$,
  since both the partial order $\leq$ and the class $\monos$ of \wmonos{}
  are defined only in terms of the \wepis{}.
\end{Rem}

\begin{Expl}
  \label{expl:Delta-Gamma-monotone}
  \begin{itemize}
  \item
    In both the \cooltechs{}
    $\Dmincooltech$ and $\Dmaxcooltech$
    on $\Delta$ defined in \autoref{sec:example-Delta}
    the partial order $\leq$
    on the objects $\numD{\nn}\in\Delta$
    is just the usual comparison of cardinalities;
    the monos are the injective maps.
    Hence $\Dmincooltech$ and $\Dmaxcooltech$
    are both \coolsmonotone{}.
  \item
    Denote by $\Gcooltech$ the \cooltech{} on $\Gamma$
    defined in \autoref{sec:cooltechs-from-factorizations}.
    It is \coolsmonotone{} since
    the \wmonos{} are opposite to the surjective maps in $\finsetp$
    and the order $\leq$
    is again just given by comparing cardinalities of finite pointed sets.
    \qedhere
  \end{itemize}
\end{Expl}

\begin{Prop}
  \label{prop:char-of-RKE-of-restriction}
  Let
  $\excooltech=(\localB,\epis,\monos)$
  be a \coolpmonotone{} \cooltech{}
  with associated quotient
  $\localAp=\localAp(\excooltech)$.
  Fix a diagram
  $\exdiagramB\colon\localB\to\extargetcat$
  in an arbitrary \inftycat{} $\extargetcat$
  and an object $\locala\in\localA$.
  \begin{enumerate}
  \item
    \label{it:inprop-partially-monotone-ll-limit}
    The functor $\exdiagramB$ is
    pointwise at $\localaB\in\localB$
    a right Kan extension
    of its restriction to $\localBll{\localaB}$
    if and only if
    \begin{equation}
      \label{eq:inprop-limit-cone-on-Es}
      \lcone{(\episneqof{\localaB})}
      \simeq\episof{\localaB}\hra\localB
      \xra{\exdiagramB}
      \extargetcat
    \end{equation}
    is a limit cone in $\extargetcat$.
  \item
    \label{it:inprop-partially-monotone-ll-vanishing}
    If the \inftycat{} $\exwaddcat\coloneqq\extargetcat$
    is weakly idempotent complete additive
    (or \pre{additive} if $\excooltech$ is \really{})
    then this happens
    if and only if
    the corresponding pointed diagram
    $\exptdiagramA\colon\localAp\to\exwaddcat$
    vanishes at $\locala$,
    \ie if and only if
    $\localevalA\exptdiagramA{\locala}$
    is a zero object in $\exwaddcat$.
  \item
    \label{it:inprop-monotone-ngeq}
    Assume that $\excooltech$ is \coolsmonotone{}.
    Then
    $\restr{\exdiagramB}{\episof{\locala}}$
    is a limit cone
    if and only if
    $\exdiagramB$ is
    pointwise at $\localaB\in\localB$
    a right Kan extension
    of its restriction to $\localBngeq{\localaB}$.
    \qedhere
  \end{enumerate}
\end{Prop}

\begin{Prf}
  We only prove
  \ref{it:inprop-partially-monotone-ll-limit}
  and
  \ref{it:inprop-partially-monotone-ll-vanishing};
  the proof of
  \ref{it:inprop-monotone-ngeq}
  is analogous to
  \ref{it:inprop-partially-monotone-ll-limit}.
  The pointwise right Kan extension of
  $\exdiagramBll{\localaB}$ at $\localaB\in\localB$
  is computed as the limit of the diagram
  \begin{equation}
    \lim
    \left(
      \undercat{\localaB}{\left(\localBll{\localaB}\right)}
      \lra
      \localB
      \xra{\exdiagramB}
      \exwaddcat
    \right)
  \end{equation}
  We show that if $\excooltech$ is \coolpmonotone{}
  then the canonical inclusion
  $
  \episneqof{\localaB}\to
  \undercat{\localaB}{\left(\localBll{\localaB}\right)}
  $
  is homotopy initial:
  \begin{itemize}
  \item
    This amounts to showing that for
    each object $\localb\in\localB$
    with $\localb<\localaB$
    and
    each arrow $\localaB\to\localb$,
    the poset of factorizations
    \begin{equation}
      \label{eq:in-proof-factorizations-through-Eneq}
      \cdtriangle
      {\localaB}{\localaB[1]}{\localb}
      {\episneq}{}{}
    \end{equation}
    is \wcontractible{}.
    The first leg in the unique $(\epis,\monos)$-factorization
    \begin{equation}
      \label{eq:in-proof-can-E-M-factorization}
      \cdtriangle
      {\localaB}{\localaB[1]}{\localb}
      {\epis}{\monos}{}
    \end{equation}
    must be \non{invertible}
    because otherwise $\locala[1]\cong\locala>\localb$
    would contradict the assumption that
    $\excooltech{}$ is \coolpmonotone{}.
    It follows that the unique factorization
    \eqref{eq:in-proof-can-E-M-factorization}
    is of type
    \eqref{eq:in-proof-factorizations-through-Eneq}
    and is therefore a terminal object in the the poset we wish to contract.
  \end{itemize}
  It follows that the desired pointwise right Kan extension
  is computed as the limit
  \begin{equation}
    \lim
    \left(
      \episneqof{\localaB}
      \lra
      \localB
      \xra{\exdiagramB}
      \exwaddcat
    \right)
  \end{equation}
  as required by
  \ref{it:inprop-partially-monotone-ll-limit}.
  Statement
  \ref{it:inprop-partially-monotone-ll-vanishing}
  now follows from
  \autoref{thm:main-thm}
  \ref{it:main-colim-lim-formula}
  which states in particular
  that the canonical map
  \begin{equation}
    \localevalB{\exdiagramB}{\localaB}
    \lra
    \lim\limits_{\localb\in\episneqof{\localaB}}
    \localevalB{\exdiagramB}{\localb}.
  \end{equation}
  is retraction with complement
  $\localevalA\exptdiagramA{\locala}$.
\end{Prf}

Fix a natural number $\kk\in\BN$.
Recall that $\Deltaleq{\kk}\subset\Delta$
denotes the full \sub{category} spanned by
the objects $\numD{\nn}$
with $\nn\leq\kk$.

\begin{dCor}
    A simplicial object
    $\exsobj\colon\Dop\to\exwaddcat$
    in a \waddjectives{} \inftycat{}
    is a left Kan extension of its restriction to
    $\Dopleq{\kk}$
    if and only if
    the corresponding connective chain complex
    $\extermcplx\in\connChof{\exwaddcat}$
    is $\kk$-truncated,
    \ie $\extermcplx_\nn\simeq 0$ for all $\nn>\kk$.
\end{dCor}

\begin{Prf}
  Apply
  \autoref{prop:char-of-RKE-of-restriction}%
  ~\ref{it:inprop-partially-monotone-ll-vanishing}
  to the \cooltech{} $\Dmincooltech$
  (or, equivalently, to the \cooltech{} $\Dmaxcooltech$)
  and dualize.
\end{Prf}

\ifoptionfinal
{}
{
  \section{Further content (not part of the final paper)}
  This section contains further material which
  will not be part of the paper itself
  but which I will include in my thesis.
  \subsection{Functoriality}
  \subfile{further-content/functoriality}

  \subsection{Computing membrane objects}
  \subfile{further-content/membrane-objects}

  \subsection{Polynomial functors}
  \subfile{further-content/polynomial-functors}
}

\bibliography{library}
\end{document}
